\documentclass[10pt,reqno,a4paper,final]{amsart}
\usepackage[british]{babel}
\usepackage{amssymb,amstext,amsthm,eucal,bbm,amscd,bbold}
\usepackage[margin=1in]{geometry}
\usepackage{array}
\usepackage{float}
\usepackage{wrapfig}
\restylefloat{figure}
\usepackage[svgnames,table]{xcolor}
\usepackage[T1]{fontenc}
\usepackage[utf8x]{inputenc}
\usepackage[small]{eulervm}
\usepackage{mathrsfs}
\usepackage{verbatim}
\usepackage{pdfsync}
\usepackage[final]{microtype}
\usepackage{booktabs} 
\usepackage{caption}
\usepackage{appendix}
\usepackage{chngcntr}
\usepackage{apptools}
\usepackage{tikz,pgfplots}
\usepgfplotslibrary{fillbetween}
\usetikzlibrary{calc,arrows,arrows.meta,cd,shapes.geometric,positioning,through,decorations.markings,decorations.text}
\tikzset{>=latex}
\usepackage[unicode,pagebackref=true]{hyperref}
\hypersetup{%
  pdftitle   = {Non-lorentzian spacetimes},
  pdfkeywords = {Newton-Cartan, Carroll, Aristotelian, G-structures,
    Klein geometry},
  pdfauthor  = {José Figueroa-O'Farrill},
  pdfcreator = {\LaTeX\ with package \flqq hyperref\frqq},
  linkcolor=NavyBlue,
  citecolor=ForestGreen,
  urlcolor=OrangeRed,
  anchorcolor=OrangeRed,
  colorlinks=true
}
\renewcommand*{\backref}[1]{}
\renewcommand*{\backrefalt}[4]{%
  \ifcase #1%
  \or [Page~#2.]%
  \else [Pages~#2.]%
  \fi%
}
\theoremstyle{plain}
\newtheorem{lemma}{Lemma}

\theoremstyle{definition}
\newtheorem{definition}[lemma]{Definition}

\newcommand{\g}{\mathfrak{g}}
\newcommand{\h}{\mathfrak{h}}
\renewcommand{\a}{\mathfrak{a}}
\renewcommand{\b}{\mathfrak{b}}
\renewcommand{\c}{\mathfrak{c}}
\renewcommand{\d}{\partial}
\newcommand{\gl}{\mathfrak{gl}}

\renewcommand{\r}{\mathfrak{r}}
\newcommand{\m}{\mathfrak{m}}
\newcommand{\n}{\mathfrak{n}}

\newcommand{\so}{\mathfrak{so}}
\newcommand{\iso}{\mathfrak{iso}}

\renewcommand{\k}{\mathfrak{k}}
\newcommand{\s}{\mathfrak{s}}

\newcommand{\be}{\boldsymbol{e}}
\newcommand{\x}{\boldsymbol{x}}
\newcommand{\y}{\boldsymbol{y}}
\newcommand{\bv}{\boldsymbol{v}}

\ifdefined\C
\renewcommand{\C}{\boldsymbol{C}}
\else
\newcommand{\C}{\boldsymbol{C}}
\fi
\renewcommand{\L}{\boldsymbol{L}}
\newcommand{\B}{\boldsymbol{B}}
\renewcommand{\P}{\boldsymbol{P}}
\newcommand{\bzero}{\boldsymbol{0}}

\newcommand{\Kgr}{\mathscr{K}}
\newcommand{\Hgr}{\mathscr{H}}

\newcommand{\eJ}{\mathscr{J}}

\newcommand{\eL}{\mathscr{L}}
\newcommand{\eX}{\mathscr{X}}
\newcommand{\eN}{\mathscr{N}}
\newcommand{\eQ}{\mathscr{Q}}

\newcommand{\Af}{\mathbb{A}}
\newcommand{\LL}{\mathbb{L}}
\newcommand{\Aff}{\operatorname{Aff}}

\newcommand{\ad}{\operatorname{ad}}
\newcommand{\Ort}{\operatorname{O}}

\newcommand{\Ad}{\operatorname{Ad}}
\newcommand{\ISO}{\operatorname{ISO}}

\newcommand{\tr}{\operatorname{tr}}

\newcommand{\RR}{\mathbb{R}}

\newcommand{\GL}{\operatorname{GL}}
\newcommand{\SO}{\operatorname{SO}}

\newcommand{\Ann}{\operatorname{Ann}}
\newcommand{\End}{\operatorname{End}}
\newcommand{\Hom}{\operatorname{Hom}}
\newcommand{\coker}{\operatorname{coker}}

\definecolor{dkgr}{rgb}{0,0.6,0}
\definecolor{gris}{rgb}{0.5,0.5,0.5}

\newcommand{\choice}[2]{\genfrac{}{}{0pt}{}{#1}{#2}}

\allowdisplaybreaks[0]
\numberwithin{equation}{section}
\AtAppendix{\counterwithin{lemma}{section}}

\begin{document}

\title{Non-lorentzian spacetimes}

\author{José Figueroa-O'Farrill}
\address{Maxwell Institute and School of Mathematics, The University
  of Edinburgh, James Clerk Maxwell Building, Peter Guthrie Tait Road,
  Edinburgh EH9 3FD, Scotland, United Kingdom}
\email{\href{mailto:j.m.figueroa@ed.ac.uk}{j.m.figueroa[at]ed.ac.uk},
  \href{https://orcid.org/0000-0002-9308-9360}{ORCID:~0000-0002-9308-9360}}
\begin{abstract}
  I review some of my recent work on non-lorentzian geometry.  I
  review the classification of kinematical Lie algebras and their
  associated Klein geometries.  I then describe the Cartan geometries
  modelled on them and their characterisation in terms of their
  intrinsic torsion.
\end{abstract}
\thanks{EMPG-21-16}
\maketitle
\tableofcontents

\section{Introduction}
\label{sec:introduction}

Ever since Minkowski introduced his eponymous spacetime (see, e.g.,
\cite{zbMATH02638586} for his famous 1908 address to the Society of
Natural Scientists in Cologne, reproduced in English translation in
\cite{MR0044931}), lorentzian geometry has played a fundamental rôle
in Physics and, until recently, has dominated our attempts to model
the universe.  Minkowski spacetime replaced the Galilei spacetime of
newtonian mechanics as the geometric arena for Einstein's theory of
special relativity, conceived to describe Maxwell's theory of
electrodynamics.  With the advent of quantum mechanics, it became the
arena of relativistic quantum mechanics and, unavoidably, of quantum
field theory.  On the other hand, Minkowski spacetime is the flat
model of lorentzian geometry, the Cartan geometry modelled on it, and
which is the basis of Einstein's theory of general relativity, which
accurately describes a wide range of gravitational phenomena.  The
attempt at marrying quantum theory and general relativity into a
quantum theory of gravity has kept the theoretical/mathematical
physics community busy for the best part of the last 75 years.  Why
then should one bother with non-lorentzian spacetimes, except as a
purely mathematical curiosity?

One answer to this question lies precisely in the difficulty to
formulate a quantum theory of gravity.  It may help to keep the
following picture in mind: the Bronstein cube of physical theories.

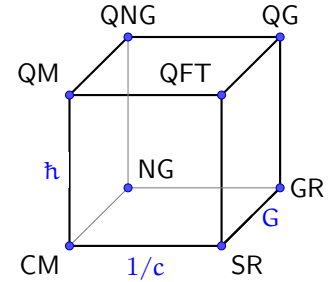
\begin{wrapfigure}{3}{0pt}
  \begin{tikzpicture}
    \coordinate [label=above:{$\mathsf{QG}$}] (qg) at (2,2,0);
    \coordinate [label=above:{$\mathsf{QNG}$}] (qng) at (0,2,0);
    \coordinate [label=above left:{$\mathsf{QM}$}] (qm) at (0,2,2);
    \coordinate [label=below left:{$\mathsf{CM}$}] (cm) at (0,0,2);
    \coordinate [label=below right:{$\mathsf{SR}$}] (sr) at (2,0,2);
    \coordinate [label=right:{$\mathsf{GR}$}] (gr) at (2,0,0);
    \coordinate [label=above left:{$\mathsf{QFT}$}] (qft) at (2,2,2);
    \coordinate [label=above right:{$\mathsf{NG}$}] (ng) at (0,0,0);
    \draw[thick] (cm) -- (qm) node [midway, left, fill=white] {\color{blue}$\hbar$};
    \draw[thick] (cm) -- (sr) node [midway, below, fill=white] {\color{blue}$1/c$};
    \draw[thick] (sr) -- (gr) node [midway, right, fill=white] {\color{blue}$G$};
    \draw[thick] (2,2,0) -- (0,2,0) -- (0,2,2) -- (2,2,2) -- (2,2,0) -- (2,0,0) -- (2,0,2);
    \draw[thick] (2,2,2) -- (2,0,2);
    \draw[gray] (2,0,0) -- (0,0,0) -- (0,2,0);
    \draw[gray] (0,0,0) -- (0,0,2);
    \foreach \point in {qg,qng,qm,cm,sr,gr,qft,ng}
    \filldraw [color=blue!70!black,fill=blue!70!white] (\point) circle (1.5pt);
  \end{tikzpicture}
  \caption{Bronstein cube}
  \label{fig:bronstein}
\end{wrapfigure}

This cube is a cartoon of physical theories one can obtain from
classical mechanics ($\mathsf{CM}$) by turning on certain parameters:
the inverse speed of light $1/c$, Newton's gravitational constant
(or, more geometrically, curvature) $G$ and Planck's constant $\hbar$.
Of course, these are physical constants and as such they take
particular values, but let us pretend that we can change them at will.
There are three directions we can go in from classical mechanics: to
special relativity ($\mathsf{SR}$) by turning on $1/c$, to quantum
mechanics ($\mathsf{QM}$) by turning on $\hbar$, and to newtonian
gravity ($\mathsf{NG}$) by turning on $G$.  From special relativity we
may go to general relativity ($\mathsf{GR}$) by turning on $G$ and to
relativistic quantum mechanics and hence quantum field theory
($\mathsf{QFT}$) by turning on $\hbar$.  These two theories are the
standard points of departure towards the final goal of a theory of
quantum gravity ($\mathsf{QG}$).  However as the picture
makes clear, there is a third possible line of approach: via the (as
yet non-existent) quantisation of newtonian gravity
($\mathsf{QNG}$).

In this review I shall remain in the non-quantum world and
concentrate on the bottom side of the Bronstein cube.  In geometrical
terms, classical mechanics and special relativity are described by
Klein geometries: Galilei and Minkowski spacetimes, whereas
turning on $G$ corresponds to constructing Cartan
geometries modelled on them.  I shall therefore start by describing
Galilei and Minkowski spacetimes and show that they are Klein
geometries of the Galilei and Poincaré groups, respectively.  I shall
recognise these groups as examples of kinematical Lie groups and
recall the classification of kinematical Lie algebras.  This will lead
to the classification of kinematical Klein geometries, of which Galilei
and Minkowski spacetimes are but two of a plethora of examples which
nevertheless give rise to a small class of Cartan geometries of
spacetimes: lorentzian, galilean (or Newton--Cartan), carrollian and
aristotelian.  I do not discuss aristotelian geometries in this
review and concentrate on the galilean and carrollian
geometries.  These Cartan geometries are examples of $G$-structures,
which I shall classify further in terms of their intrinsic torsion.

This short review is organised as follows.  In
Section~\ref{sec:two-models-universe} I discuss two classical models
of the universe, as motivation for the rest of the geometries in the
review: Galilei spacetime is discussed in
Section~\ref{sec:galilei-spacetime} and Minkowski spacetime in
Section~\ref{sec:minkowski-spacetime}.  Comparison of the relativity
groups of Galilei and Minkowski spacetimes, suggests the notion of a
kinematical Lie group and in Section~\ref{sec:klas} I define the
notion of a kinematical Lie algebra and review their classification,
given in Table~\ref{tab:KLAs}.  This is a necessary first step to the
determination of the associated kinematical Klein geometries.  These
are discussed in Section~\ref{sec:kinem-klein-geom} and their
classification (up to coverings) is given in
Table~\ref{tab:spacetimes}.  After a brief review of the basic notions
of Klein geometry in Section~\ref{sec:basic-notions-klein}, I discuss
the four classes of kinematical Klein geometries in the Table:
lorentzian in Section~\ref{sec:lorentz-klein-geom}, riemannian in
Section~\ref{sec:riem-klein-geom}, galilean in
Section~\ref{sec:galil-klein-geom} and carrollian in
Section~\ref{sec:carr-klein-geom}.  I do not dwell on the riemannian
case, since riemannian manifolds do not play the rôle of spacetimes.
For the lorentzian and carrollian cases, I give explicit geometric
realisations of the spacetimes in Table~\ref{tab:spacetimes},
postponing the galilean case until later in the review.
In Section~\ref{sec:kinem-cart-geom} I turn my attention to
the kinematical Cartan geometries modelled on the kinematical Klein
geometries.  In Section~\ref{sec:basic-notions-cartan} I present a
very brief review of the basic notions of Cartan geometry that I
shall require, mainly about $G$-structures since the kinematical Cartan
geometries turn out to be $G$-structures.  I also review the notion
of the intrinsic torsion of a $G$-structure, which will allow a
refinement of the classification of kinematical Cartan geometries.
Section~\ref{sec:lorentzian-geometry} briefly recaps, mostly for
orientation, the case of lorentzian geometry.  In
Section~\ref{sec:newt-cart-geom} I discuss Newton--Cartan geometry,
which is the Cartan geometry modelled on galilean spacetimes and I
summarise the classification in terms of intrinsic torsion.  In Section
\ref{sec:null-reductions} I discuss how to obtain Newton--Cartan
geometries as null reductions of lorentzian geometries and this is
used in Section~\ref{sec:another-look-at} in order to give the
promised geometric realisations of the galilean Klein geometries in
Table~\ref{tab:spacetimes}.  In
Section~\ref{sec:carrollian-geometry} I discuss carrollian geometry
and summarise the classification in terms of intrinsic torsion.  The
nomenclature of carrollian structures is meant to be reminiscent of
the classification of hypersurfaces in riemannian geometry.  This is
made precise in Section~\ref{sec:null-hypersurfaces}, where I show
that a null hypersurface in a lorentzian manifold has a carrollian
structure and re-interpret the intrinsic torsion in this light using
the properties of the null Weingarten map and associated null second
fundamental form.  In Section~\ref{sec:bundles-scales} I discuss
another source of carrollian structures: namely, the natural
structure on the bundle of scales of a riemannian conformal
manifold.  Finally, in Section~\ref{sec:conclusion-outlook} I give a
brief overview of related topics not covered in this review and some
natural extensions of the work described here.

\section*{Acknowledgments}

I would like to thank Dmitri Alekseevsky for the invitation to
contribute this review.  I would also like to thank the following
past, present and hopefully future collaborators for many hours of
stimulating discussions on non-lorentzian and, more generally, Cartan
geometry: Tomasz Andrzejewski, Andrew Beckett, Eric Bergshoeff,
Joaquim Gomis, Can Görmez, Ross Grassie, Jelle Hartong, Emil Have, Yannick
Herfray, Stefan Prohazka, Jakob Salzer, Andrea Santi, Dennis The and
Dieter Van den Bleeken.

\section{Two models of the universe}
\label{sec:two-models-universe}

In this section I review two classical models of the universe:
Galilei spacetime of newtonian mechanics and Minkowski spacetime of
special relativity.  Both spacetimes are described by an affine space,
homogeneous under the action of a kinematical Lie group (to be defined
below), but their invariant structures differ: whereas Galilei
spacetime has a Galilei-invariant clock and ruler, Minkowski
spacetime has a Poincaré-invariant proper distance.  In the
corresponding Cartan geometries, the clock and ruler will be seen as a
flat example of a (weak) Newton--Cartan structure, whereas the proper
distance is, of course, a flat example of a lorentzian metric.

Although the physically relevant dimension is $4$, let us work
in $d+1$ dimensions.  Let $\Af^{d+1}$ denote ($d+1$)-dimensional real
affine space.  In the present context, points in $\Af^{d+1}$ are
referred to as \textbf{(spacetime) events}.  The affine space
$\Af^{d+1}$ is a torsor over $\RR^{d+1}$, thought of as an abelian
group under vector addition.  The action of $\RR^{d+1}$ on $\Af^{d+1}$
is via parallel displacements.  Given any two points
$a,b \in \Af^{d+1}$, there exists a unique parallel displacement
$v \in \RR^{d+1}$ such that $b = a + v$.  It is customary to refer to
$v$ as $b-a$ and hence to identify parallel displacements with
differences of points.

It is often convenient for calculations to use an explicit model for
$\Af^{d+1}$ as the affine hyperplane in $\RR^{d+2}$ with equation
$x^{d+2}=1$.  This embeds the affine group $\Aff(d+1,\RR)$ into the
general linear group $\GL(d+2,\RR)$ as the subgroup which preserves
that hyperplane or, more explicitly, as the subgroup consisting of
matrices of the form
\begin{equation}\label{eq:affine-trans}
  \begin{pmatrix}
    L & v \\ 0 & 1
  \end{pmatrix}
\end{equation}
where $v \in \RR^{d+1}$ and $L \in \GL(d+1,\RR)$.  The relativity
groups of Galilei and Minkowski spacetimes are subgroups of the affine
group consisting of all the parallel displacements but a restricted
subgroup of linear transformations consisting of orthogonal
transformations and boosts; although of course the notion of boost
differs in the galilean and lorentzian worlds.

\subsection{Galilei spacetime}
\label{sec:galilei-spacetime}

The following formulation of the Galilei spacetime is originally due
to Weyl \cite{MR988402}.  Galilei spacetime is a triple $(\Af^{d+1},
\tau, \lambda)$ where the \textbf{clock} $\tau : \RR^{d+1} \to \RR$
and the \textbf{ruler} $\lambda : \ker \tau \to \RR$ are defined as
follows.

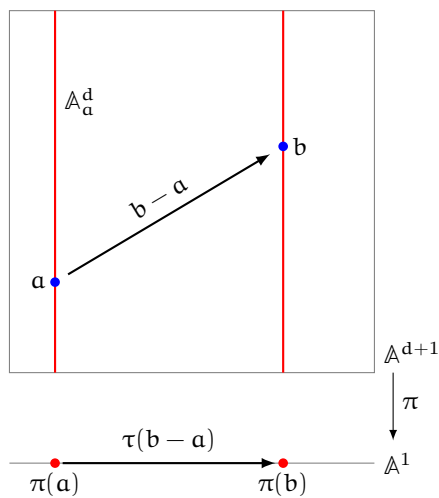
\begin{figure}[h!]
  \centering
    \begin{tikzpicture}[x=1.0cm,y=1.0cm,scale=0.6]
      \coordinate [label=above right:{$\Af^{d+1}$}] (a4) at (4,1); 
      \coordinate [label=right:{$\Af^1$}] (a1) at (4,-1);
      \coordinate [label=below:{$\pi(a)$}] (pa) at (-3,-1); 
      \coordinate [label=below:{$\pi(b)$}] (pb) at (2,-1); 
      \coordinate [label=left:{$a$}] (a) at (-3,3); 
      \coordinate [label=right:{$b$}] (b) at (2,6);
      \coordinate [label=right:{$\Af^d_a$}] (a3) at (-3,7);
      \draw [black!50!white, thin] (-4,1) -- (4,1) -- (4,9) -- (-4,9) -- cycle;
      \draw [black!50!white, thin] (-4,-1) -- (4,-1);
      \draw [red, thick] (-3,1) -- (-3,9);
      \draw [red, thick] (2,1) -- (2,9);
      \draw[->,shorten >=2mm,shorten <=2mm,black,thick] (a)--(b) node[midway,sloped,above]{$b-a$}; 
      \draw[->,shorten >=1mm,shorten <=1mm,black,thick] (pa)--(pb) node[midway,above]{$\tau(b-a)$}; 
      \begin{scope}[transform canvas={xshift=0.7em}]
        \draw [->, shorten >=3mm,thin, black] (a4) -- node[above right] {$\pi$} (a1);
      \end{scope}
      \foreach \point in {pa,pb}
      \fill [red] (\point) circle (3pt);
      \foreach \point in {a,b}
      \fill [blue] (\point) circle (3pt);
  \end{tikzpicture}
  \caption{The clock fibration $\pi : \Af^{d+1} \to \Af^1$}
  \label{fig:clock-fibration}
\end{figure}

The clock measures the time interval $\tau(b-a)$ between two events
$a,b \in \Af^{d+1}$.  In the explicit model
$\Af^{d+1} \subset \RR^{d+2}$, letting $a = (\x, x^{d+1}, 1)$ and
$b = (\y, y^{d+1},1)$, with $\x,\y \in \RR^d$, then $\tau(b-a) =
y^{d+1} - x^{d+1}$.  Two events $a,b \in \Af^{d+1}$ are
\textbf{simultaneous} if $\tau(b-a) = 0$.  Fixing an event $a$, the
set of events simultaneous to $a$ defines a $d$-dimensional affine
subspace
\begin{equation}
  \Af^d_a := a + \ker \tau = \left\{a + v ~ \middle | ~\tau(v) =
    0\right\}
\end{equation}
of $\Af^{d+1}$.  The quotient $\Af^{d+1}/\ker \tau$ is an affine
line $\Af^1$, so that the clock gives a fibration $\pi: \Af^{d+1}
\to \Af^1$ whose fibre at $\pi(a)$ consists of the affine
hypersurface $\Af^d_a$, as illustrated in
Figure~\ref{fig:clock-fibration}.

The ruler measures the (euclidean) distance between simultaneous
events.  If $a,b \in \Af^{d+1}$ are simultaneous, $\lambda(b-a)$ is
the euclidean length of $b - a \in \ker\tau$.  Again, in the explicit
model, if $a = (\x, x^{d+1}, 1)$ and $b = (\y, y^{d+1},1)$, with
$\x,\y \in \RR^d$ and $x^{d+1}=y^{d+1}$, are any two simultaneous
events, then $\lambda(b-a) = \|\y - \x\|$, the euclidean distance
between $\x$ and $\y$.

The relativity group of Galilei spacetime is called the
\textbf{Galilei group} and it consists of those affine transformations
of $\Af^{d+1}$ which preserve the clock and the ruler.  It embeds in
$\GL(d+2,\RR)$ as those matrices of the form
\begin{equation}
  \label{eq:gal-in-gl}
  \begin{pmatrix}
    R & \bv & \boldsymbol{p} \\
    0 & 1 & s \\
    0 & 0 & 1
  \end{pmatrix},
\end{equation}
where $R \in \Ort(d)$, $\boldsymbol{p},\bv \in \RR^d$ and $s \in
\RR$.

The action of the matrix in equation~\eqref{eq:gal-in-gl} on an event
$(\x, t, 1)$ gives the event $(R\x + t \bv + \boldsymbol{p}, t +
s,1)$ which may be interpreted as the composition of an orthogonal
transformation $\x \mapsto R \x$, a \textbf{Galilei boost} $\x
\mapsto \x + t \bv$, a spatial translation $\x \mapsto \x +
\boldsymbol{p}$ and a temporal translation $t \mapsto t + s$.  The
stabiliser subgroup of the event $(\bzero,0,1)$, isomorphic to what is
commonly called the \textbf{homogeneous Galilei group}, consists of
orthogonal transformations and boosts.  As discussed later, the Cartan
geometry modelled on Galilei spacetime is an $G$-structure with
$G$ the homogeneous Galilei group.

The Lie algebra of the Galilei group is unsurprisingly called the
\textbf{Galilei algebra} and it is isomorphic to the subalgebra of
$\gl(d+2,\RR)$ consisting of matrices of the form
\begin{equation}
  \label{eq:gal-algebra-gl}
  \begin{pmatrix}
    A & \bv & \boldsymbol{p}\\
    0 & 0 & s\\
    0 & 0 & 0
  \end{pmatrix},
\end{equation}
where $A \in \so(d)$, $\bv,\boldsymbol{p} \in \RR^d$ and $s \in \RR$.
Introduce a basis $L_{ab} = - L_{ba}, B_a, P_a, H$ by
\begin{equation}
  \label{eq:gal-basis}
  \begin{pmatrix}
    A & \bv & \boldsymbol{p}\\
    0 & 0 & s\\
    0 & 0 & 0
  \end{pmatrix} = \tfrac12 A^{ab} L_{ab} + v^a B_a + p^a P_a + s H
\end{equation}
and in this way easily work out the Lie brackets of the Galilei
algebra in this basis.  The nonzero brackets are given by
\begin{equation}
  \label{eq:gal-algebra-brackets}
  \begin{split}
    [L_{ab},L_{cd}] &= \delta_{bc} L_{ad} - \delta_{ac} L_{bd} -  \delta_{bd} L_{ac} + \delta_{bd} L_{ac} \\
    [L_{ab}, B_b] &= \delta_{bc} B_a - \delta_{ac} B_b\\
    [L_{ab}, P_b] &= \delta_{bc} P_a - \delta_{ac} P_b\\
    [B_a, H] &= P_a.
  \end{split}
\end{equation}
This shows that $L_{ab}$ span an $\so(d)$ subalgebra, relative to
which $B_a,P_a$ transform according to the three-dimensional vector
representation and $H$ transforms as the one-dimensional scalar
representation.  These turn out to be the defining properties of
a kinematical Lie algebra (with spatial isotropy).

\subsection{Minkowski spacetime}
\label{sec:minkowski-spacetime}

Minkowski spacetime is described by a pair $(\Af^{d+1}, \Delta)$, where
the \textbf{proper distance} $\Delta : \RR^{d+1} \to
\RR$ is defined as follows.  In the explicit model of $\Af^{d+1}
\subset \RR^{d+2}$, if $a = (\x,t,1)$ and $b = (\y,s,1)$,
\begin{equation}
  \label{eq:proper-distance}
  \Delta(b-a) = \|\y - \x\|^2 - c^2 (s-t)^2,
\end{equation}
where $c$ is a parameter interpretable as the speed of light.

There is no longer a separate clock and ruler, or in Minkowski's own
words (in the English translation of \cite{MR0044931}):
\begin{quotation}
  Henceforth space by itself, and time by itself, are doomed to fade
  away into mere shadows, and only a kind of union of the two will
  preserve an independent reality.
\end{quotation}
In particular, there is no longer an invariant notion of simultaneity
between events.  Simultaneity is replaced by a notion of causality,
geometrised by lightcones at every spacetime event $a$: the \textbf{lightcone}
$\LL_a$ of $a$ being defined as those events which are a zero proper distance
away from $a$:
\begin{equation}
  \LL_a = \left\{ b \in \Af^{d+1} ~\middle | ~\Delta(b-a) = 0\right\}.
\end{equation}
Two events $a,b \in \Af^{d+1}$ are said to be causally related if
$\Delta(b-a) \leq 0$.  If $a = (\x, t, 1)$ and $b = (\y, s, 1)$ are
causally related, one says that $a$ is in the causal future of $b$ (and
hence $b$ is in the causal past of $a$) if $t - s > 0$.

The relativity group of Minkowski spacetime is the \textbf{Poincaré
  group} and consists of those affine transformations which preserve
the proper distance between events.  It embeds in $\GL(d+2,\RR)$ as
those matrices
\begin{equation}
  \label{eq:poin-in-gl}
  \begin{pmatrix}
    L & v \\
    0 & 1
  \end{pmatrix}
\end{equation}
where $L^T \eta L = \eta$ and $v \in \RR^{d+1}$.  The Poincaré group
is thus isomorphic to the semidirect product $\Ort(d,1) \ltimes
\RR^{d+1}$, where $\Ort(d,1)$ is the \textbf{Lorentz group}.  The
effect of the Poincaré transformation with matrix
\eqref{eq:poin-in-gl} on an event $(x,1)$ is the event $(Lx + v, 1)$,
which is the effect of a Lorentz transformation $x \mapsto L x$ and a
(spatiotemporal) translation $x \mapsto x + v$.  The Lorentz group is
the stabiliser subgroup of the event $(\bzero,0,1)$ and, of course,
lorentzian geometry is the study of $\Hgr$-structures with $\Hgr$ the
Lorentz group.

The Lie algebra of the Poincaré group embeds in $\gl(d+2,\RR)$ as those
matrices of the form
\begin{equation}
  \begin{pmatrix}
    A & v \\
    0 & 0
  \end{pmatrix}
\end{equation}
where $A^T \eta + \eta A = 0$ and $v \in \RR^{d+1}$.  Introducing a basis
$L_{mn} = - L_{nm}, P_m$, where now $m,n = 0,1,\dots, d$, by
\begin{equation}
  \begin{pmatrix}
    A & v \\
    0 & 0
  \end{pmatrix} = \tfrac12 A^{mn} L_{mn} + v^m P_m,
\end{equation}
it is easy to calculate the nonzero Lie brackets:
\begin{equation}\label{eq:poincare-brackets}
  \begin{split}
    [L_{mn},L_{pq}] &= \eta_{np} L_{mq} - \eta_{mp} L_{nq} - \eta_{nq} L_{mp} + \eta_{mq} L_{np}\\
    [L_{mn},P_p] &= \eta_{np} P_m - \eta_{mp} P_n.
  \end{split}
\end{equation}
To ease comparison with the Galilei algebra
\eqref{eq:gal-algebra-brackets}, let $P_m = (P_a, H = P_0)$ and
$L_{mn} = (L_{ab}, B_a = L_{0a})$, relative to which the brackets
become
\begin{equation}
  \label{eq:poincare-kla-brackets}
  \begin{split}
    [L_{ab},L_{cd}] &= \delta_{bc} L_{ad} - \delta_{ac} L_{bd} -  \delta_{bd} L_{ac} + \delta_{bd} L_{ac} \\
    [L_{ab}, B_b] &= \delta_{bc} B_a - \delta_{ac} B_b\\
    [L_{ab}, P_b] &= \delta_{bc} P_a - \delta_{ac} P_b\\
    [B_a, B_b] &= c^2 L_{ab}\\
    [B_a, P_b] &= \delta_{ab} H\\
    [B_a, H] &= c^2 P_a.
  \end{split}
\end{equation}
Again $L_{ab}$ span an $\so(d)$ subalgebra relative to which $B_a,P_a$
transform according to the vector representation and $H$ transforms
according to the one-dimensional scalar representation.  What
distinguishes the Poincaré and Galilei algebras are the Lie brackets
which do not involve the $L_{ab}$: the last bracket in
equation~\eqref{eq:gal-algebra-brackets} and the last three brackets
in equation~\eqref{eq:poincare-kla-brackets}.

The Lie brackets in \eqref{eq:poincare-kla-brackets} depend explicitly
on the parameter $c$, the speed of light, which may formally be set to
any desired value.  For any nonzero value, the resulting Lie algebras
are isomorphic: simply rescale $B_a \mapsto c^{-1}B_a$ and $H
\mapsto c^{-1}H$, which is an isomorphism for nonzero $c$, resulting in
the brackets with $c=1$.  On the other hand, setting $c=0$ results in
a non-isomorphic Lie algebra, with brackets
\begin{equation}
  \label{eq:carroll-kla-brackets}
  \begin{split}
    [L_{ab},L_{cd}] &= \delta_{bc} L_{ad} - \delta_{ac} L_{bd} -  \delta_{bd} L_{ac} + \delta_{bd} L_{ac} \\
    [L_{ab}, B_b] &= \delta_{bc} B_a - \delta_{ac} B_b\\
    [L_{ab}, P_b] &= \delta_{bc} P_a - \delta_{ac} P_b\\
    [B_a, P_b] &= \delta_{ab} H.
  \end{split}
\end{equation}
This algebra was first studied by Lévy-Leblond \cite{MR0192900}, who
named it the \textbf{Carroll algebra} in honour of Lewis Carroll, the
pseudonym of Charles Dodgson, the author of \emph{Alice's Adventures
  in Wonderland}.

Alternatively, one can obtain the Galilei algebra formally as the
limit $c \to \infty$ of the Poincaré algebra; that is, as a Lie
algebra contraction.  Indeed, let us rescale $B_a \mapsto c^{-2}B_a$,
relative to which the Lie brackets become
\begin{equation}
  \label{eq:poincare-pre-contraction}
  \begin{split}
    [L_{ab},L_{cd}] &= \delta_{bc} L_{ad} - \delta_{ac} L_{bd} -  \delta_{bd} L_{ac} + \delta_{bd} L_{ac} \\
    [L_{ab}, B_b] &= \delta_{bc} B_a - \delta_{ac} B_b\\
    [L_{ab}, P_b] &= \delta_{bc} P_a - \delta_{ac} P_b\\
    [B_a, B_b] &= c^{-2} L_{ab}\\
    [B_a, P_b] &= c^{-2}\delta_{ab} H\\
    [B_a, H] &= P_a.
  \end{split}
\end{equation}
The rescaling is an isomorphism for any non-zero value of $c^{-2}$ and
hence results in an isomorphic Lie algebra.  In the limit $c^{-2} \to
0$, the rescaling is singular, but the Lie brackets do have a limit
and it is evident that the resulting Lie brackets in this limit agree
with those of the Galilei algebra in
equation~\eqref{eq:gal-algebra-brackets}, showing that the
Galilei algebra is a contraction of the Poincaré algebra.

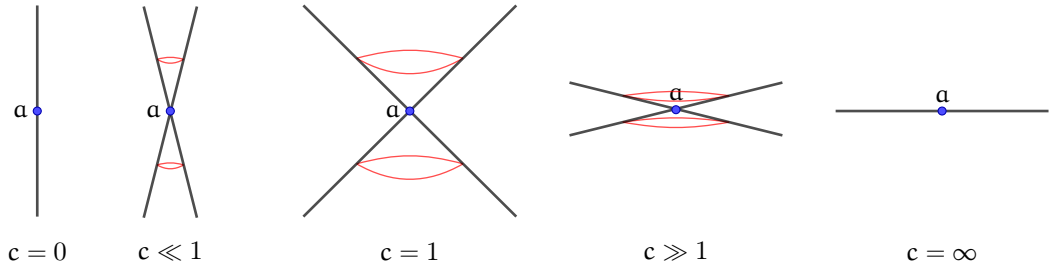
\begin{figure}[h!]
  \centering
  \begin{tikzpicture}[x=0.7cm,y=0.7cm]
    %
    %
    %
    %
    \coordinate (sw) at (6.,-1.);
    \coordinate (se) at (8.,-1.);
    \coordinate (ne) at (8.,1.);
    \coordinate (nw) at (6.,1.);
    \coordinate (usw) at (2.2503423410911236,-1.0203486648922608);
    \coordinate (use) at (2.7503423410911236,-1.0203486648922608); 
    \coordinate (une) at (2.7503423410911236,0.9796513351077392); 
    \coordinate (unw) at (2.2503423410911236,0.9796513351077392); 
    \coordinate (nsw) at (10.999892248241657,-0.21078008816633959);
    \coordinate (nse) at (12.999892248241657,-0.21078008816633959); 
    \coordinate (nne) at (12.999892248241657,0.28921991183366036);
    \coordinate (nnw) at (10.999892248241657,0.28921991183366036);
    \coordinate [label=left:{$a$}] (c) at (0,0); 
    \coordinate [label=left:{$a$}] (u) at (2.5,0); 
    \coordinate [label=left:{$a$}] (l) at (7,0); 
    \coordinate [label=above:{$a$}] (n) at (12,0.025); 
    \coordinate [label=above:{$a$}] (g) at (17,0); 
    %
    %
    \begin{scope}[line width=0.5pt,opacity=0.7,color=red]
      \draw (nw) to [out=330,in=210] (ne);
      \draw (sw) to [out=330,in=210] (se); 
      \draw (nw) to [out=15,in=165] (ne);      
      \draw (sw) to [out=15,in=165] (se);      
      \draw (unw) to [out=330,in=210] (une);
      \draw (usw) to [out=330,in=210] (use); 
      \draw (unw) to [out=15,in=165] (une);      
      \draw (usw) to [out=15,in=165] (use);      
      \draw (nnw) to [out=350,in=190] (nne);
      \draw (nsw) to [out=350,in=190] (nse);
      \draw (nnw) to [out=8,in=172] (nne);      
      \draw (nsw) to [out=8,in=172] (nse);      
    \end{scope}
    \begin{scope}[line width=1pt,opacity=0.7,color=black]
      \draw (5.,2.)-- (9.,-2.);
      \draw (9.,2.)-- (5.,-2.);
      \draw (9.999892248241657,0.5392199118336602)-- (13.999892248241657,-0.46078008816633964);
      \draw (13.999892248241657,0.5392199118336602) -- (9.999892248241657,-0.46078008816633964);
      \draw (3.0003423410911236,1.979651335107739)-- (2.0003423410911236,-2.020348664892261);
      \draw (2.0003423410911236,1.979651335107739)-- (3.0003423410911236,-2.020348664892261);
      \draw (15.,0.)-- (19.,0.);
      \draw (0.,2.)-- (0.,-2.);
    \end{scope}
    %
    %
    \node at (7,-3) [above] {$c=1$};
    \node at (2.5,-3) [above] {$c\ll 1$};
    \node at (0,-3) [above] {$c=0$};
    \node at (12,-3) [above] {$c\gg 1$};
    \node at (17,-3) [above] {$c=\infty$};
    %
    %
    \foreach \point in {c,u,l,n,g}
    \filldraw [color=blue!70!black,fill=blue!70!white] (\point) circle (1.5pt);
  \end{tikzpicture}
  \caption{Cartoon of the effect of varying $c$ on the lightcone $\LL_a$}
  \label{fig:lightcones}
\end{figure}

These two limits of the Poincaré algebra: $c \to 0$ and $c \to
\infty$ can be understood geometrically according to what they do to
the lightcones.  In the limit $c \to \infty$, the lightcone $\LL_a$ 
opens up to become the affine hyperplane of events simultaneous to
$a$, where now $\Af^{d+1}$ is to be interpreted as Galilei spacetime.
In the limit $c \to 0$, the lightcone $\LL_a$ closes up to become the
affine temporal line based at $a$, as depicted in
Figure~\ref{fig:lightcones}.  Taking the limit $c \to \infty$ says
that any characteristic speed in the physics is much smaller than the
speed of light and hence this limit is typically known as the
non-relativistic limit.  In contrast, in the limit $c \to 0$, since no
material body can travel faster than the speed light, motion is
impossible.  This is called the ultra-relativistic (or ultra-local)
limit.

\section{Kinematical Lie algebras}
\label{sec:klas}

In Section~\ref{sec:two-models-universe} I discussed two models of
the universe: the Galilei and Minkowski spacetimes.  Both are affine
spaces homogeneous under the action of a ``kinematical'' Lie group:
the Galilei and Poincaré groups, respectively.  In this section I
will define this notion and discuss the classification of kinematical
Lie algebras.

More than half a century ago, Bacry and Lévy-Leblond \cite{MR0238545}
asked themselves the question of which were the
possible kinematics, rephrasing the question mathematically as the
classification of kinematical Lie algebras.  A careful comparison of
the Poincaré, Galilei and Carroll algebras encountered in
Section~\ref{sec:two-models-universe} suggests the following
definition, now for ($d+1$)-dimensional spacetimes.\footnote{Strictly
  speaking, the definition is for spatially isotropic spacetimes.
  There are generalisations where the rotational subalgebra $\r$ in
  the definition is replaced by a Lorentz subalgebra $\so(d-1,1)$ or
  more generally a pseudo-orthogonal subalgebra $\so(p,d-p)$.  Such
  homogeneous spaces do occur in nature.  Indeed, as shown in
  \cite{Gibbons:2019zfs}, the blow-up \cite{Ashtekar:1978zz} of
  spatial infinity of Minkowski spacetime is a homogeneous space of
  the Poincaré group with lorentzian isotropy.  There are other
  homogeneous spaces of the Poincaré group occurring at the asymptotic
  infinities of Minkowski spacetime, as discussed in
  \cite{Figueroa-OFarrill:2021sxz}, which have carrollian or
  Carroll-like structures.} 

\begin{definition}\label{def:KLA}
  Let $V$ be a $d$-dimensional euclidean vector space and $\so(V)$ the
  corresponding orthogonal Lie algebra.  A \textbf{kinematical Lie
    algebra} (for spatially isotropic ($d+1$)-dimensional spacetimes)
  is a real Lie algebra $\k$ with a subalgebra $\r \cong \so(V)$ and
  such that under the restriction to $\r$ of the adjoint
  representation of $\k$,
  \begin{equation}
    \k \cong \so(V) \oplus 2V \oplus \RR,
  \end{equation}
  where $\so(V)$, $V$ and $\RR$ are the adjoint, defining and trivial
  one-dimensional representations of $\so(V)$, respectively.
\end{definition}

In addition, Bacry and Lévy-Leblond also imposed that $\k$ should
admit two automorphisms: \emph{parity} $P_a \mapsto -P_a$ and
\emph{time-reversal} $H \mapsto -H$; although they did point out that
those restrictions were ``by no means compelling'' and indeed twenty
years later, Bacry and Nuyts \cite{MR857383} relaxed them to arrive at
the classification of four-dimensional kinematical Lie algebras as we
understand them today.  This classification was recovered using
deformation theory in \cite{Figueroa-OFarrill:2017ycu} and extended to
arbitrary dimension in
\cite{Figueroa-OFarrill:2017tcy,Andrzejewski:2018gmz}.  The case of
$d=1$ is a recontextualisation of the Bianchi classification of
three-dimensional real Lie algebras \cite{Bianchi,MR1900159}, here
re-interpreted as kinematical Lie algebras for two-dimensional
spacetimes.  The cases of $d=2$ and $d=3$ are complicated by the
existence of additional $\r$-invariant tensors in $\wedge^2 V$ and
$\wedge^3 V$, respectively, which can contribute to the brackets.
And, indeed, there are kinematical Lie algebras in dimension $2+1$ and
$3+1$ which have no higher-dimensional analogues.  I will refer the
interested reader to the papers cited above and will concentrate here
on those kinematical Lie algebras which exist in generic dimensions.

It is often convenient to choose an orthonormal basis for $V$ and
corresponding generators $L_{ab}=-L_{ba}, B_a, P_a, H$, with
$a,b=1,2,\dots, d$, for $\k$ in terms of which, the defining
properties of a kinematical Lie algebra are contained in the following
brackets:
\begin{equation}
  \label{eq:kla-generic}
  \begin{split}
    [L_{ab},L_{cd}] &= \delta_{bc} L_{ad} - \delta_{ac} L_{bd} -  \delta_{bd} L_{ac} + \delta_{bd} L_{ac}\\
    [L_{ab}, B_b] &= \delta_{bc} B_a - \delta_{ac} B_b\\
    [L_{ab}, P_b] &= \delta_{bc} P_a - \delta_{ac} P_b\\
    [L_{ab}, H] &= 0.
  \end{split}
\end{equation}

A naive (but effective) way to approach the classification of
kinematical Lie algebras is simply to write down the most general
$\r$-invariant Lie brackets for the generators $B_a,P_a,H$ and impose
the Jacobi identity.  The Jacobi identity cuts out an algebraic
variety $\eJ$ in $\Hom_{\r}(\wedge^2 W, \k)$, for $W = 2 V \oplus
\RR$.  Two points in $\eJ$ define isomorphic kinematical Lie algebras
if and only if they are in the same orbit of $\GL_{\r}(W)$, the group
of $\r$-invariant general linear transformation of $W$.  One studies
the orbit decomposition of $\eJ$ and selects a unique representative
for each orbit.

This procedure results in Table~\ref{tab:KLAs}, which lists the
kinematical Lie algebras in generic dimension $d+1$.  For $d \leq 2$,
there are some degeneracies (e.g., if $d=2$, the Galilei algebra $\g$
is isomorphic to the Carroll algebra $\c$), but for general $d$ the
table below lists non-isomorphic kinematical Lie algebras and for
$d>3$ the table is complete.  The table lists the nonzero Lie brackets
in addition to the defining ones in equation~\eqref{eq:kla-generic}.
It also uses a shorthand notation omitting indices.  The only
$\r$-invariant tensor which can appear is $\delta_{ab}$ and hence
there is an unambiguous way to add indices.  There is no standard
notation for all the kinematical Lie algebras, so I have made some
choices.

\begin{table}[h!]
  \centering
  \caption{Kinematical Lie algebras in generic dimension}
  \label{tab:KLAs}
  \rowcolors{2}{blue!10}{white}
  \begin{tabular}{>{$}l<{$}|*{5}{>{$}l<{$}}|l} \toprule
    \multicolumn{1}{c|}{Name} &  \multicolumn{5}{c|}{Nonzero Lie brackets in addition to those in \eqref{eq:kla-generic}} & \multicolumn{1}{c}{Comments} \\\midrule
       \s & & & & & & \\
       \g & [H,\B] = - \P & & & & & \\
    \n^0    & [H,\B] = \B + \P & [H, \P] = \P & & & & \\
    \n^+_\gamma  & [H,\B] = \gamma \B & [H,\P] = \P & & & & $\gamma \in [-1,1]$ \\
    \n^-_\chi & [H,\B] = \chi \B + \P & [H,\P] = \chi \P - \B &  & & & $\chi \geq 0$ \\
    \c & & & & [\B,\P] = H & & \\
    \choice{\iso(d,1)}{\iso(d+1)} & [H,\B] = -\varepsilon \P & &  [\B,\B]= \varepsilon \L & [\B,\P] = H & & $\varepsilon = \pm 1$ \\
    \so(d+1,1) & [H,\B] = \B & [H,\P] = -\P & &  [\B,\P] = H + \L & & \\
    \choice{\so(d,2)}{\so(d+2)} & [H,\B] = -\varepsilon \P & [H,\P] = \varepsilon \B &  [\B,\B]= \varepsilon \L & [\B,\P] = H &  [\P,\P] = \varepsilon \L & $\varepsilon = \pm 1$ \\ \bottomrule
  \end{tabular}
\end{table}

Let us now describe each of the algebras in turn:
\begin{itemize}
\item The Lie algebra $\s$ is the \textbf{static} kinematical Lie
  algebra: all additional brackets are zero.  Therefore every
  kinematical Lie algebras is a deformation of $\s$.

\item The Galilei algebra is denoted by $\g$ and a closely related Lie
  algebra has been denoted by $\n^0$.  In $\g$ and $\n^0$, $\ad_H$
  is not semisimple, but has a nontrivial Jordan block:
  \begin{equation}
    \ad_H^{\g}
    \begin{pmatrix} \B \\ \P \end{pmatrix} = \begin{pmatrix}
      0 & -1 \\ 0 & 0
    \end{pmatrix} \begin{pmatrix} \B \\ \P
    \end{pmatrix}\qquad\text{and}\qquad
    \ad_H^{\n^0}
    \begin{pmatrix} \B \\ \P \end{pmatrix} = \begin{pmatrix}
      1 & 1 \\ 0 & 1
    \end{pmatrix} \begin{pmatrix} \B \\ \P
    \end{pmatrix}.  
  \end{equation}
  
\item There are two one-parameter families of Lie algebras, deforming
  the \textbf{Newton--Hooke} algebras \cite{MR0238545,MR334711}:
  \begin{itemize}
  \item $\n^+_\gamma$,  with $\gamma \in [-1,1]$, which for $\gamma =
    -1$ is the Newton--Hooke algebra $\n^+$ in the notation of
    \cite{MR334711}; and
  \item $\n^-_\chi$, with $\chi \geq 0$, which for $\chi = 0$ is the
    other Newton--Hooke algebra $\n^-$, in the notation of \cite{MR334711}.
  \end{itemize}
  These two families correspond to the cases where
  $\ad_H$ is semisimple, with real eigenvalues in $\n^+_\gamma$ and
  complex eigenvalues in $\n^-_\chi$:
  \begin{equation}
  \ad_H^{\n^+}
  \begin{pmatrix} \B \\ \P \end{pmatrix} = \begin{pmatrix}
    \gamma & 0 \\ 0 & 1
  \end{pmatrix} \begin{pmatrix} \B \\ \P
  \end{pmatrix}\qquad\text{and}\qquad
  \ad_H^{\n^-}
  \begin{pmatrix} \B \\ \P \end{pmatrix} = \begin{pmatrix}
    \chi & 1 \\ -1 & \chi
  \end{pmatrix} \begin{pmatrix} \B \\ \P
  \end{pmatrix}.  
\end{equation}

\item The Carroll algebra \cite{MR0192900} is denoted $\c$.

\item The Poincaré algebra is $\iso(d,1)$ and the euclidean algebra is
  $\iso(d+1)$.

\item The remaining Lie algebras are $\so(d+2)$, $\so(d+1,1)$ and
  $\so(d,2)$, which for $d \geq 2$ are semisimple.  Finite-dimensional
  semisimple Lie algebras are rigid, but they do admit contractions.
  Of the Lie algebras in the Table, those which can be obtained as
  contractions from these are $\s$, $\c$, $\g$, $\n^-_{\chi =0}$,
  $\n^+_{\gamma =-1}$, $\iso(d+1)$ and $\iso(d,1)$, which are
  precisely the kinematical Lie algebras admitting parity and
  time-reversal automorphisms and whose four-dimensional avatars were
  the objects in the first classification of Bacry and Lévy-Leblond
  \cite{MR0238545}.
\end{itemize}

\section{Kinematical Klein geometries}
\label{sec:kinem-klein-geom}

In this section I review the classification of kinematical Klein
geometries arrived at in \cite{Figueroa-OFarrill:2018ilb} and further
studied in \cite{Figueroa-OFarrill:2019sex}.

\begin{definition}
  A \textbf{kinematical Lie group} is a real Lie group whose Lie
  algebra is kinematical (see Definition~\ref{def:KLA}).  A
  \textbf{kinematical Klein geometry} is a ($d+1$)-dimensional
  homogeneous space of such a kinematical Lie group, whose Klein pair
  $(\k,\h)$ is such that $\k$ is a kinematical Lie algebra and $\h$ is
  a subalgebra containing the subalgebra $\r \cong \so(V)$ and such
  that under the restriction to $\r$ of the adjoint representation,
  $\h \cong \r \oplus V$.
\end{definition}

Already in the pioneering work of Bacry and Lévy-Leblond
\cite{MR0238545}, for $d=3$ and the restricted list of
kinematical Lie algebras admitting parity and time-reversal
automorphisms, a list of eleven possible kinematical Klein pairs are
discussed.  Upon closer analysis one of their Klein pairs is not
effective and describes a static aristotelian spacetime, resulting in
ten kinematical Klein geometries: all of which are reductive and
symmetric. The classification of (simply-connected) kinematical Klein
geometries was arrived at in \cite{Figueroa-OFarrill:2018ilb}, to
where I refer the interested reader for details.  The classification
for $d\geq 3$ is summarised in Table~\ref{tab:spacetimes}, where the
Lie brackets of the kinematical Lie algebra $\k$ are listed in a basis
where the subalgebra $\h$ is spanned by $L_{ab},B_a$.

\begin{table}[h!]
  \centering
  \caption{($d+1$)-dimensional kinematical Klein geometries ($d\geq 3$)}
  \label{tab:spacetimes}
  \rowcolors{2}{blue!10}{white}
  \resizebox{\textwidth}{!}{
    \begin{tabular}{l|>{$}l<{$}|*{5}{>{$}l<{$}}}\toprule
      \multicolumn{1}{c|}{Name} & \multicolumn{1}{c|}{Klein pair} & \multicolumn{5}{c}{Nonzero Lie brackets in addition to \eqref{eq:kla-generic}} \\\midrule
      Minkowski & (\iso(d,1),\so(d,1)) & [H,\B] = -\P & & [\B,\B] = \L & [\B,\P] = H &\\
      de~Sitter & (\so(d+1,1),\so(d,1)) & [H,\B] = -\P & [H,\P] = -\B & [\B,\B]= \L & [\B,\P] = H & [\P,\P]= - \L \\
      anti~de~Sitter & (\so(d,2),\so(d,1))  & [H,\B] = -\P & [H,\P] = \B & [\B,\B]= \L & [\B,\P] = H & [\P,\P] = \L \\\midrule
      euclidean & (\iso(d+1),\so(d+1)) &[H,\B] = \P & & [\B,\B] = -\L & [\B,\P] = H &  \\
      sphere &  (\so(d+2),\so(d+1)) & [H,\B] = \P & [H,\P] = -\B & [\B,\B]= -\L & [\B,\P] = H & [\P,\P]= - \L  \\
      hyperbolic &  (\so(d+1,1),\so(d+1)) & [H,\B] = \P & [H,\P] = \B & [\B,\B]= -\L & [\B,\P] = H & [\P,\P] = \L \\\midrule
      Galilei & (\g,\iso(d)) & [H,\B] = -\P & & & & \\
      de~Sitter--Galilei & (\n^+_{\gamma =-1},\iso(d))  & [H,\B] = -\P & [H,\P] = -\B & & & \\
      torsional de~Sitter--Galilei & (\n^+_{\gamma\in(-1,1)},\iso(d)) & [H,\B] = -\P & [H,\P] = \gamma\B + (1+\gamma)\P & & & \\
      torsional de~Sitter--Galilei & (\n^0,\iso(d)) & [H,\B] = -\P & [H,\P] = \B + 2\P & & & \\
      anti~de~Sitter--Galilei & (\n^-_{\chi=0},\iso(d)) & [H,\B] =  -\P & [H,\P] = \B & & & \\
      torsional anti~de~Sitter--Galilei & (\n^-_{\chi>0},\iso(d)) & [H,\B] = -\P & [H,\P] = (1+\chi^2) \B + 2\chi \P & & & \\\midrule
      Carroll  & (\c,\iso(d)) & & & & [\B,\P] = H & \\
      de~Sitter--Carroll & (\iso(d+1), \iso(d)) & & [H,\P] = -\B & & [\B,\P] = H & [\P,\P] = -\L \\
      anti~de~Sitter--Carroll & (\iso(d,1), \iso(d)) & & [H,\P] = \B & & [\B,\P] = H & [\P,\P] = \L  \\
      lightcone & (\so(d+1,1), \iso(d)) & [H,\B] = \B & [H,\P] = -\P & & [\B,\P] = H + \L & \\\bottomrule
    \end{tabular}
  }
\end{table}

The table is divided into four sections, from top to bottom:
lorentzian, riemannian, galilean and carrollian Klein geometries.  In
order to explain this coarser classification I need to explain these
terms, but before doing so, I will review briefly some basic notions
of Klein geometry.

\subsection{Basic notions of Klein geometry}
\label{sec:basic-notions-klein}

Let $\k$ be a Lie algebra and $\h$ a Lie subalgebra.  The Klein pair
$(\k,\h)$ is said to be
\begin{itemize}
\item \textbf{effective} if $\h$ does not contain a nonzero ideal of
  $\k$, and
\item \textbf{geometrically realisable} if there exists a kinematical
  Lie group $\Kgr$ with Lie algebra $\k$ such that the connected
  subgroup $\Hgr$ generated by $\h$ is closed.  The homogeneous space
  $\Kgr/\Hgr$ is a \textbf{geometric realisation} of $(\k,\h)$.
\end{itemize}

Two Klein pairs $(\k_1,\h_1)$ and $(\k_2,\h_2)$ are
\textbf{isomorphic} if there is a Lie algebra isomorphism
$\varphi: \k_1 \to \k_2$ with $\varphi(\h_1) = \h_2$.  The fundamental
theorem of Klein geometry states that isomorphism classes of
geometrically realisable, effective Klein pairs are in bijective
correspondence with isomorphism classes of simply-connected
homogeneous spaces.  Paraphrasing, effective and geometrically
realisable Klein pairs classify homogeneous spaces up to coverings.

Given an effective, geometrically realisable Klein pair $(\k,\h)$, the
\textbf{linear isotropy representation} $\lambda: \h \to \gl(\k/\h)$
is defined as the representation induced by the restriction to $\h$ of
the adjoint representation of $\k$.  Explicitly, letting $X \in \h$
and $\overline Y \in \k/\h$ denote the residue class modulo $\h$ of $Y
\in \k$,
\begin{equation}
  \lambda_X \overline Y := \overline{\ad_XY}.
\end{equation}

The holonomy principle establishes a bijective correspondence between
invariant tensors of the linear isotropy representation and tensor
fields on the homogeneous space which are invariant under the
kinematical Lie group $\Kgr$.  It is these invariant tensor fields
which determine the type of the Klein geometry.

A Klein pair $(\k,\h)$ is said to be \textbf{reductive} if the short
exact sequence
\begin{equation}
  \begin{tikzcd}
    0 \arrow[r] & \h \arrow[r] & \k \arrow[r] & \k/\h \arrow[r] & 0
  \end{tikzcd}
\end{equation}
splits in the category of $\h$-modules.  This is equivalently to the
existence of a complement $\m$ to $\h$ in $\k$ which is stable under
the restriction to $\h$ of the adjoint action of $\k$; that is,
$\k = \h \oplus \m$ with $[\h,\m] \subset \m$ in the obvious notation.
A reductive Klein pair $(\k = \h \oplus \m, \h)$ is said to be
\textbf{symmetric} if $[\m,\m] \subset \h$.  All Klein geometries in
Table~\ref{tab:spacetimes} are reductive with the exception of the
lightcone and of the reductive ones, all but the ``torsional'' ones
are symmetric.

The ``torsional'' adjective refers to the torsion of the canonical
invariant connection on a reductive Klein geometry.  An
\textbf{invariant connection} on a Klein geometry with Klein pair
$(\k,\h)$ is an $\h$-equivariant linear map $\Lambda: \k \to
\gl(\k/\h)$, written $X \mapsto \Lambda_X$, whose restriction to $\h$
is the linear isotropy representation; that is, for all $X \in \h$,
$\Lambda_X = \lambda_X$.

All Klein geometries in Table~\ref{tab:spacetimes} admit invariant
connections, with the exception of the lightcone (for $d\geq 2$).
These have been tabulated in \cite{Figueroa-OFarrill:2019sex}.

If $(\k = \h \oplus \m,\h)$ is reductive, then invariant connections
are determined by their \textbf{Nomizu map} \cite{MR0059050}, the
$\h$-equivariant bilinear map $\alpha : \m \times \m \to \m$ defined
by $\alpha(X,Y) := \Lambda_X Y$ for all $X,Y \in \m$.  The torsion
$\Theta : \wedge^2\m \to \m$ and curvature $\Omega: \wedge^2\m \to
\gl(\m)$ of an invariant connection are given in terms of the Nomizu
map by
\begin{equation}
  \begin{split}
    \Theta(X,Y) &= \alpha(X,Y) - \alpha(Y,X) - [X,Y]_\m\\
    \Omega(X,Y)Z &= \alpha(X,\alpha(Y,Z)) - \alpha(Y,\alpha(X,Z)) -
    \alpha([X,Y]_\m,Z) - [[X,Y]_\h, Z],
  \end{split}
\end{equation}
where $X = X_\h + X_\m$ is the decomposition of $X \in \k = \h \oplus
\m$.

The \textbf{canonical invariant connection} of a reductive Klein
geometry is the unique invariant connection with zero Nomizu map.  Its
torsion and curvature take particularly simple forms:
\begin{equation}
  \Theta(X,Y) = - [X,Y]_\m \qquad\text{and}\qquad \Omega(X,Y)Z = -
  [[X,Y]_\h, Z].
\end{equation}
It follows that a reductive Klein geometry is torsion-free if and only
if it is symmetric.

The holonomy representation of the canonical invariant connection is
isomorphic to the ideal $[\m,\m]_\h \subset \h$ acting on $\m$ via the
restriction of the linear isotropy representation
\cite[§12]{MR0059050}.  It follows that any $\Kgr$-invariant tensor
field on the homogeneous space is parallel with respect to to the
canonical invariant connection.

\subsection{Lorentzian Klein geometries}
\label{sec:lorentz-klein-geom}

The lorentzian Klein geometries all have Klein pairs $(\k,\h)$ with
$\h \cong \so(d,1)$.  They can be characterised by the existence of an
$\h$-invariant lorentzian inner product on $\k/\h$; that is, a
nondegenerate bilinear form in $\left( \odot^2 (\k/\h)^* \right)^\h$
with lorentzian signature.  The holonomy principle gives rise to a
$\Kgr$-invariant lorentzian metric on the homogeneous space.  Since
$\dim \Kgr = \tfrac12 (d+1)(d+2)$, the homogeneous space has constant
sectional curvature, so one of Minkowski spacetime (flat), de~Sitter
spacetime (positive curvature) and anti~de~Sitter spacetime (negative
curvature).  Geometrically, there is a one-parameter family of (anti)
de~Sitter spacetimes, corresponding to the value of the scalar
curvature, but all de Sitter spacetimes are isomorphic as homogeneous
spaces and similarly for all anti~de~Sitter spacetimes.

Let us give a geometric realisation for each one.  Minkowski
spacetime, as described in Section~\ref{sec:minkowski-spacetime}, is a
geometric realisation of the Klein pair $(\iso(d,1), \so(d,1)$.

A geometric realisation of ($d+1$)-dimensional de~Sitter spacetime is
provided by a quadric hypersurface in a lorentzian vector space of
dimension $d+2$.  Let $x^0,x^1,\dots,x^{d+1}$ be Cartesian coordinates
for $\RR^{d+2}$ with the flat lorentzian metric
\begin{equation}\label{eq:flat-lor-metric}
  -(dx^0)^2 + \sum_{i=1}^{d+1} (dx^i)^2.
\end{equation}
Then for any $\ell \neq 0$, the hypersurface cut out by the quadric
\begin{equation}
  -(x^0)^2 + \sum_{i=1}^{d+1} (dx^i)^2 = \ell^2
\end{equation}
is a geometric realisation of the Klein pair $(\so(d+1,1),\so(d,1))$.

A geometric realisation of ($d+1$)-dimensional anti~de~Sitter
spacetime is also provided by a quadric hypersurface, but this time in
a pseudo-euclidean vector space with an inner product of signature
$(d,2)$.  Let $x^0,x^1,\dots,x^{d+1}$ be Cartesian coordinates
for $\RR^{d+2}$ with the flat metric
\begin{equation}
  - (dx^0)^2 + \sum_{i=1}^d (dx^i)^2 - (dx^{d+1})^2.
\end{equation}
Then for all $\ell \neq  0$, the hypersurface cut out by the quadric
\begin{equation}
  - (x^0)^2 + \sum_{i=1}^d (x^i)^2 - (x^{d+1})^2 = -\ell^2
\end{equation}
is a geometric realisation of the Klein pair $(\so(d,2),\so(d,1))$.

\subsection{Riemannian Klein geometries}
\label{sec:riem-klein-geom}

The riemannian Klein geometries all have Klein pairs $(\k,\h)$ with
$\h \cong \so(d+1)$.  They can be characterised by the existence of an
invariant euclidean inner product on $\k/\h$.  The story is very
similar to the lorentzian case above.  The resulting homogeneous
spaces are riemannian symmetric spaces with constant sectional
curvature: euclidean space (flat), the round sphere (positive
curvature) and hyperbolic space (negative curvature).  Although
strictly speaking they are Klein geometries of kinematical Lie groups,
they play no rôle as spacetimes.  One could eliminate them from the
discussion by imposing a further condition on a kinematical Klein
geometry; namely, that the generic orbits of the one-parameter
subgroup corresponding to any $v^a B_a \in \h$ on the homogeneous
space should be non-compact.  In the riemannian Klein geometries,
these one-parameter subgroups act as rotations and therefore have
compact orbits.

\subsection{Galilean Klein geometries}
\label{sec:galil-klein-geom}

The galilean Klein geometries all have Klein pairs $(\k,\h)$ with
$\h \cong \iso(d)$.  They can be characterised by the existence of an
invariant covector in the dual $(\k/\h)^*$ of the linear isotropy
representation and an invariant symmetric bivector in
$\odot^2(\k/\h)$.  The corresponding invariant tensor fields on the
homogeneous space can be interpreted as a clock and ruler, just as in
Galilei spacetime in Section~\ref{sec:galilei-spacetime}.  All
galilean Klein geometries are reductive and therefore admit a
canonical invariant connection. This connection is torsion-free for
the Galilei and (anti) de~Sitter--Galilei spacetimes, which explains
the adjective ``torsional'' in the other two classes.

Galilei spacetime, as described in
Section~\ref{sec:galilei-spacetime}, is a geometric realisation of the
Klein pair $(\g, \iso(d)$.  I will postpone discussion of the
geometric realisations of the other galilean Klein geometries until
Section~\ref{sec:null-reductions} after I discuss how to to obtain
Newton--Cartan structures by null reductions of lorentzian manifolds.

\subsection{Carrollian Klein geometries}
\label{sec:carr-klein-geom}

The carrollian Klein geometries are in a certain sense dual to the galilean
Klein geometries.  They all have Klein pairs $(\k,\h)$ where $\h \cong
\iso(d)$, but the $\iso(d)$ subalgebra of $\gl(\k/\h)$ in the
carrollian and galilean cases are not conjugate under $\GL(\k/\h)$.
Indeed, they have different invariant tensors of the linear isotropy
representation.  Carrollian Klein geometries are characterised by the
existence of an invariant vector in $\k/h$ and an invariant symmetric
bilinear form in $\odot^2(\k/\h)^*$.  Except for the lightcone, which
is not reductive, all other carrollian Klein geometries are reductive
and symmetric.  Anti~de~Sitter--Carroll is a Klein geometry associated
to the Poincaré algebra and is the geometry underlying the blow-up of
past and future timelike infinities in Minkowski spacetime
\cite{Figueroa-OFarrill:2021sxz}.

All of the carrollian Klein geometries may be realised geometrically
as null hypersurfaces in lorentzian manifolds.  As the name suggests,
the ($d+1$)-dimensional lightcone is realised geometrically as the
future (or past) deleted lightcone in a ($d+2$)-dimensional lorentzian
vector space.  Let $x^0,x^1,\dots,x^{d+1}$ be cartesian coordinates
for $\RR^{d+2}$ and consider the flat lorentzian metric given in
equation~\eqref{eq:flat-lor-metric}.  Then the null hypersurface
\begin{equation}
    -(x^0)^2 + \sum_{i=1}^{d+1} (dx^i)^2 = 0
\end{equation}
with $x^0 > 0$, say, is a geometric realisation of the Klein pair
$(\so(d+1,1),\iso(d))$.

The other three carrollian Klein geometries may be realised as null
hypersurfaces in Minkowski spacetime (for the Carroll spacetime), de
Sitter spacetime (for the de~Sitter--Carroll spacetime) and
anti~de~Sitter spacetime (for the anti~de~Sitter--Carroll spacetime).

Consider ($d+2$)-dimensional Minkowski spacetime with coordinates
$x^0,x^1,\dots,x^{d+1}$ and metric
\begin{equation}
  -(dx^0)^2 + \sum_{i=1}^{d+1} (dx^i)^2.
\end{equation}
As shown originally in \cite{Duval:2014uoa}, Carroll spacetime embeds
here as the null hyperplane $x^0 = x^{d+1}$.  Indeed, it is not hard
to show (see, e.g., \cite[§4.2.5]{Figueroa-OFarrill:2018ilb}) that
this null hyperplane is a geometric realisation of the Klein pair
$(\c,\iso(d))$.

Consider now ($d+3$)-dimensional Minkowski spacetime with coordinates
$x^0,x^1,\dots,x^{d+2}$ and metric
\begin{equation}
  -(dx^0)^2 + \sum_{i=1}^{d+2} (dx^i)^2.
\end{equation}
Let $\eQ$ denote the quadric hypersurface defined by
\begin{equation}
  -(x^0)^2 + \sum_{i=1}^{d+2} (x^i)^2 = \ell^2,
\end{equation}
which, as seen in Section~\ref{sec:lorentz-klein-geom} and for every
$\ell \neq 0$, is covered by de~Sitter spacetime.  Let $\eN$ denote the 
hyperplane cut out by $x^0 = x^{d+2}$.  Then as shown in
\cite[§4.2.5]{Figueroa-OFarrill:2018ilb}, the intersection
$\eQ \cap \eN$ is a geometric realisation of the Klein pair
$(\iso(d+1),\iso(d))$.

Consider finally a ($d+3$)-dimensional pseudo-euclidean space with
coordinates $x^0,x^1,\dots,x^{d+2}$ and metric
\begin{equation}
  -(dx^0)^2 + \sum_{i=1}^{d+1} (dx^i)^2 - (dx^{d+2})^2.
\end{equation}
Let $\eQ'$ denote the quadric hypersurface defined by
\begin{equation}
  -(x^0)^2 + \sum_{i=1}^{d+1} (x^i)^2 - (x^{d+2})^2 = -\ell^2,
\end{equation}
which, as seen in Section~\ref{sec:lorentz-klein-geom} and for any
$\ell \neq 0$, is covered by anti~de~Sitter spacetime.  Let $\eN'$
denote the hyperplane cut out by $x^{d+1}=x^{d+2}$.  Then as shown in 
\cite[§4.2.5]{Figueroa-OFarrill:2018ilb} and \cite{Morand:2018tke},
the intersection $\eQ' \cap \eN'$ is a geometric realisation of the
Klein pair $(\iso(d,1),\iso(d))$.  This geometric realisation has
been used recently in order to describe the asymptotic geometry of
Minkowski spacetime in terms of homogeneous spaces of the Poincaré
group \cite{Figueroa-OFarrill:2021sxz}.

\section{Kinematical Cartan geometries}
\label{sec:kinem-cart-geom}

In this section I discuss the Cartan geometries modelled on the
kinematical Klein geometries discussed in
Section~\ref{sec:kinem-klein-geom}, but before doing so, I review the
basic notions of Cartan geometry relevant to our discussion.  A good
treatment of Cartan geometry is given in \cite{MR1453120} and
there is a growing list of explicit applications of Cartan geometry to
gravitation
\cite{Wise:2006sm,Herfray:2020rvq,Herfray:2021xyp,Herfray:2021qmp}.

\subsection{Basic notions of Cartan geometry}
\label{sec:basic-notions-cartan}

To specify a Cartan geometry, one needs, in addition to an effective
Klein pair $(\k,\h)$, a Lie group $\Hgr$ with Lie algebra $\h$ and an
action of $\Hgr$ on $\k$ extending the adjoint action on $\h$ and
denoted $\Ad$.  More precisely, a \textbf{Cartan geometry} modelled on
$(\k,\h)$ with group $\Hgr$ is a right $\Hgr$-principal bundle $P \to
M$ together with a \textbf{Cartan connection}: a $\k$-valued one-form
$\omega$ on $P$ satisfying the following conditions:
\begin{enumerate}
\item (\emph{non-degeneracy}) for each $p \in P$, $\omega_p : T_pP \to
  \g$ is an isomorphism;
\item (\emph{equivariance}) for every $h \in \Hgr$, $R^*_h \omega =
  \Ad(h^{-1}) \circ \omega$, where $R_h$ is the diffeomorphism of $P$
  induced by the right action of $h \in \Hgr$; and
\item (\emph{normalisation}) for all $X \in \h$, $\omega(\xi_X) = X$,
  where $\xi_X \in \eX(P)$ is the corresponding fundamental vector
  field.
\end{enumerate}

The induced action of $\Hgr$ on $\k/\h$ is the linear isotropy
representation of $\Hgr$.  If faithful, the Cartan geometry is said to
be of the \textbf{first order} and it then follows that the principal
bundle $P \to M$ is a reduction of the frame bundle of $M$; that is,
an $G$-structure with $G=\Hgr$.  This will always be the case for the
kinematical Cartan geometries under discussion.

$G$-structures come with additional structure: a soldering form.  In
the present context it is given by $\theta \in \Omega^1(P;\k/\h)$, the
projection to $\k/\h$ of the Cartan connection.  The soldering form is
both horizontal and equivariant and hence it defines a section of
$\Hom(TM,P\times_{\Hgr} \k/\h)$.  The non-degeneracy condition of the
Cartan connection says that the soldering form defines a bundle
isomorphism between the tangent bundle $TM$ and the ``fake tangent
bundle'' $P\times_{\Hgr} \k/\h$.  This allows to identify tensor
bundles over $M$ with the corresponding associated vector bundles to
the principal bundle $P\to M$, a fact that shall be used tacitly.

Let $\omega\in\Omega^1(P,\h)$ be an Ehresmann connection on $P$.  If
the Cartan geometry is reductive in addition to being of the first
order, then the $\h$-component of the Cartan connection is such an
Ehresmann connection, but Ehresmann connections exist even if the
Cartan geometry is not reductive.  An Ehresmann connection
defines a Koszul connection on any associated vector bundle to
$P \to M$ and, in particular, on the fake tangent bundle.  The
soldering form transports that Koszul connection to an affine
connection on $TM$.  Affine connections obtained in this way are said
to be \textbf{adapted} to the $G$-structure.

The difference between two adapted connections is a one-form with
values in the adjoint bundle $\Ad P = P \times_{\Hgr} \h$.  The linear
isotropy representation $\lambda : \h \to \gl(\k/\h)$ allows us to
view $\Ad P$ as a sub-bundle of the endomorphisms of the fake tangent
bundle and, via the soldering form, as endomorphisms of the tangent
bundle.  In summary, whereas the difference between any two affine
connections on $TM$ is a one-form with values in $\End(TM)$, if the
affine connections are adapted, it is a one-form taking values in the
sub-bundle $\Ad P$.

Let us introduce the notation $V := \k/\h$ and let us define the
\textbf{Spencer differential}
\begin{equation}
  \d : \Hom(V,\h) \to \Hom(\wedge^2V, V) \qquad\text{by}\qquad (\d
  \kappa)(v,w) = \kappa_v w - \kappa_w v
\end{equation}
for all $v,w \in V$.  As with every linear map, $\d$ fits inside a
four-term exact sequence
\begin{equation}
  \begin{tikzcd}
    0 \arrow[r] & \ker \d \arrow[r] & \h \otimes V^* \arrow[r,"\d"] &
    V \otimes \wedge^2 V^* \arrow[r] & \coker \d \arrow[r] & 0
  \end{tikzcd}
\end{equation}
and hence an exact sequence of the corresponding associated
vector bundles on $M$.  The sections of the associated vector bundles
to these representations can be interpreted as follows:
\begin{itemize}
\item $V \otimes \wedge^2 V^*$: torsion of (adapted) affine connections;
\item $\h \otimes V^*$: difference between adapted affine connections;
\item $\ker \d$: differences which do not alter the torsion; and
\item $\coker \d$: \textbf{intrinsic torsions} of adapted connections,
  so called because this is part of the torsion which does not depend
  on the connection and hence it is intrinsic to the $G$-structure.
\end{itemize}

The intrinsic torsion is a coarse yet easy to determine invariant of
the kinematical Cartan geometries.  A more detailed treatment can be
found in \cite{Figueroa-OFarrill:2020gpr}.

\subsection{Lorentzian geometry}
\label{sec:lorentzian-geometry}

Let $V = \RR^{d+1}$ thought of as the defining representation of
$\Ort(d,1)$, which is the subgroup of $\GL(d+1,\RR)$ leaving a
lorentzian inner product invariant.  If $M$ is a ($d+1$)-dimensional
manifold, a $\Ort(d,1)$-structure is a sub-bundle of the frame bundle
consisting of the pseudo-orthonormal frames.  Every such frame gives
isomorphism $T_p M \to \RR^{d+1}$ for each point $p$ where it is
defined and the lorentzian inner product on
$\RR^{d+1}$ pulls back to a lorentzian inner product on $T_p M$.
Restricting to pseudo-orthonormal frames, these inner products define a
smooth lorentzian metric $g$ on $M$.  The fundamental theorem of
riemannian geometry, the existence of a unique torsion-free metric
connection, may be re-interpreted in the language of $G$-structures
as saying that any adapted (i.e., metric) connection on the
orthonormal frame bundle can be modified to have zero torsion (since
$\coker \d = 0$) and, moreover, that there is a unique such
modification (since $\ker\d = 0$).

\subsection{Newton--Cartan geometry}
\label{sec:newt-cart-geom}

Galilean $G$-structures and their adapted connections were first
discussed in \cite{MR175340, MR334831} and further studied in
\cite{Bernal:2002ph, Bekaert:2014bwa, Bekaert:2015xua,
  Figueroa-OFarrill:2020gpr}.

Let $V = \RR^{d+1}$.  It is convenient to choose a suggestive notation
for the elementary basis of $\RR^{d+1}$: $\be_0,\be_1,\dots,\be_d$ and
for the canonical dual basis $\alpha^0,\alpha^1,\dots,\alpha^d$.  Let
$G$ denote the subgroup of $\GL(V)$ which fixes $\alpha^0 \in
V^*$ and $\sum_{a=1}^d \be_a \be_a \in \odot^2 V$.  Explicitly,
\begin{equation}
  \label{eq:G-gal}
  G = \left\{
    \begin{pmatrix}1 & \bzero^T \\ \bv & A\end{pmatrix} ~ \middle | ~  \bv
  \in \RR^d,~A \in O(d)\right\} < \GL(d+1,\RR),
\end{equation}
with Lie algebra
\begin{equation}
  \label{eq:g-gal}
  \g = \left\{
    \begin{pmatrix}0 & \bzero^T \\ \bv & A\end{pmatrix} ~ \middle | ~  \bv
  \in \RR^d,~A \in \so(d)\right\} < \gl(d+1,\RR).
\end{equation}

Let $M$ be an ($d+1$)-dimensional manifold with a $G$-structure with
$G$ the group in equation~\eqref{eq:G-gal}.  The characteristic tensor
fields are now a nowhere-vanishing one-form $\tau \in \Omega^1(M)$
(the \textbf{clock one-form}) and a corank-one positive-semidefinite
symmetric bivector $\lambda \in \Gamma(\odot^2 TM)$ (the
\textbf{spatial cometric}) with $\lambda(\tau,-) = 0$.  The triple
$(M,\tau,\lambda)$ is a \textbf{(weak) Newton--Cartan geometry}, which
can be promoted to a \textbf{Newton--Cartan geometry} by the addition
of an adapted affine connection $\nabla$.

Already in \cite{MR334831}, the kernel and cokernel of the Spencer
differential for a Newton--Cartan $G$-structure was determined to be
isomorphic to $\wedge^2 V^*$ as a $G$-module and that under the
isomorphism $\coker \d \cong \wedge^2 V^*$, the intrinsic torsion of
an adapted connection is sent to $d\tau \in \Omega^2(M)$.  It was
further shown in \cite{Figueroa-OFarrill:2020gpr} that in generic $d$
(here, $d\neq 1,4$) there are three $G$-submodules of $\coker \d$ and,
correspondingly, three kinds of Newton--Cartan geometries, in the
notation of \cite[Table~I]{Christensen:2013lma}, who first identified
these classes:
\begin{enumerate}
\item \textbf{torsionless Newton--Cartan geometry} (NC) if $d\tau = 0$;
\item \textbf{twistless torsional Newton--Cartan geometry} (TTNC), if
  $d\tau \wedge \tau = 0$; and
\item \textbf{torsional Newton--Cartan geometry} (TNC), if $d\tau
  \wedge \tau \neq 0$.
\end{enumerate}
If $d=1$, $d\tau \wedge \tau = 0$ on dimensional grounds, hence there
are two intrinsic torsion classes of two-dimensional Newton--Cartan
geometries: those where $d\tau = 0$ and those where $d\tau \neq 0$.
Similarly, if $d=4$, and if the $G$-structure reduces further to a
$G_0$-structure, with $G_0 \subset G$ the connected component of the
identity -- for example, if $M$ were simply-connected -- then the
subbundle $\ker \tau \subset TM$ is oriented and riemannian, so that
one can distinguish between those five-dimensional torsional
Newton--Cartan geometries where the restriction of $d\tau$ to $\ker
\tau$ is self-dual, anti~self-dual or neither.

\subsubsection{Null reductions}
\label{sec:null-reductions}

In this section I review the construction of (weak) Newton--Cartan
geometries via null reductions of lorentzian geometries
\cite{PhysRevD.31.1841,Julia:1994bs}.

Let $(N,g)$ be a lorentzian manifold with a nowhere-vanishing null
Killing vector field $\xi \in \eX(N)$.  Let us assume that $\xi$ is
complete and that it generates a one-parameter subgroup $\Gamma$ of
isometries of $(N,g)$ in such a way that the quotient $M = N/\Gamma$
is a smooth manifold.  Let $\pi : N \to M$ denote the resulting
principal $\Gamma$ bundle.  As is well known, the pullback
$\pi^* : \Omega^\bullet(M) \to \Omega^\bullet(N)$ of differential
forms sets up a $C^\infty(M)$-module isomorphism between
$\Omega^\bullet(M)$ and the \textbf{basic forms}:
\begin{equation}
  \label{eq:basic-forms}
  \Omega^\bullet_\Gamma(N) = \left\{ \omega \in \Omega^\bullet(N)
    ~\middle | ~ \imath_\xi \omega = 0\quad\text{and}\quad \imath_\xi
    d\omega = 0\right\}.
\end{equation}
The first condition ($\imath_\xi \omega = 0$) says that $\omega$ is
\textbf{horizontal} and the second condition, together with
horizontality, says that $\eL_\xi \omega = 0$, so that $\omega$ is
$\Gamma$-invariant.  Let $\flat: TN \to T^*N$ and $\sharp: T^*N \to
TN$ denote the musical isomorphisms induced by the lorentzian metric
$g$.  I shall use the same notation for the corresponding
$C^\infty(M)$-module isomorphisms between the spaces of sections:
$\flat : \eX(N) \to \Omega^1(N)$ and $\sharp: \Omega^1(N) \to
\eX(N)$.

The one-form $\xi^\flat \in \Omega^1(N)$ metrically dual to the null
Killing vector field is basic: it is horizontal because $\xi$ is null
and it is invariant because $\xi$ is Killing.  Therefore $\xi^\flat =
\pi^*\tau$ for some nowhere vanishing $\tau \in \Omega^1(M)$.

Let $\alpha,\beta \in \Omega^1(M)$ and let $X=(\pi^*\alpha)^\sharp,
Y=(\pi^*\beta)^\sharp \in \eX(N)$.  Let $f = g(X,Y) \in C^\infty(N)$.
Since $\xi$ is Killing, it follows that $\xi(f) = 0$, so that $f =
\pi^*h$ for some $h \in C^\infty(M)$.  Define $\lambda \in
\Gamma(\odot^2 M)$ by $\lambda(\alpha,\beta) = h$.  In other words,
$\lambda$ is defined by
\begin{equation}
  \pi^* \lambda(\alpha,\beta) = g((\pi^*\alpha)^\sharp,(\pi^*\beta)^\sharp).
\end{equation}
It is easy to see that $\lambda$ is well-defined, obeys $\lambda(\tau,-)
= 0$, is positive-semidefinite and has corank $1$.  In other words,
$(M,\tau,\lambda)$ is a (weak) Newton--Cartan geometry.

\subsubsection{Another look at the galilean Klein geometries}
\label{sec:another-look-at}

As promised, I now describe geometric realisations of the galilean
Klein geometries in Section~\ref{sec:galil-klein-geom} as null
reductions of lorentzian manifolds.  There results are joint with
Stefan Prohazka and Ross Grassie and will appear in a forthcoming
paper on Bargmann spacetimes.  In \cite{Figueroa-OFarrill:2018ilb},
their geometric realisability was shown non-constructively.

Let us consider Klein pairs $(\b, \h)$, where $\b$ is a \textbf{generalised
Bargmann algebra}, namely a one-dimensional extension (not necessarily
central) of a kinematical Lie algebra $\k$:
\begin{equation}
  \begin{tikzcd}
    0 \arrow[r] & \RR Z \arrow[r] & \b \arrow[r] & \k \arrow[r] & 0,
  \end{tikzcd}
\end{equation}
with $Z$ the additional generator.  The name of the Lie algebras
is due to the fact that the Bargmann algebra is the universal central
extension of the Galilei algebra, with additional bracket:
\begin{equation}
  \label{eq:bargmann-generic}
  [B_a, P_b] = \delta_{ab} Z.
\end{equation}
Every generalised Bargmann algebra has a basis
$L_{ab}, B_a, P_a, H, Z$ and shares the kinematical Lie
brackets~\eqref{eq:kla-generic} in addition to
\eqref{eq:bargmann-generic}.  It follows by $\r$-equivariance that $Z$
transforms under the trivial one-dimensional representation of
$\so(V)$.  The Newton--Hooke algebras also admit
central extensions, whereas the kinematical Lie algebras $\n^0$,
$\n^+_\gamma$, for $\gamma \in (-1,1]$ and $n^-_\xi$ for $\xi>0$ admit
non-central extensions.

Let us define the generalised Bargmann algebras $\b^+_\gamma$, for
$\gamma \in [-1,1)$, $\b^-_\chi$, for $\chi \geq 0$ and $\b^0$ as in
Table~\ref{tab:gen-bargmann}, which lists all nonzero Lie brackets in
addition to those in \eqref{eq:kla-generic} and
\eqref{eq:bargmann-generic}.  The first three Lie algebras in the
Table are the universal central extensions of the Galilei algebra $\g$
and the Newton--Hooke algebras $\n^\pm$, whereas the last three are
non-central extensions of $\n^+_{\gamma\in(-1,1)}$, $\n^0$ and
$\n^-_{\chi>0}$, respectively. Notice that even when $\b^0$ agrees
with $\b^+_{\gamma = 1}$, it is an extension of $\n^0$ and not of
$\n^+_{\gamma =1}$.

\begin{table}
  \centering
  \caption{Some relevant generalised Bargmann algebras}
  \label{tab:gen-bargmann}
  \begin{tabular}{>{$}l<{$}|*{3}{>{$}l<{$}}}
    \multicolumn{1}{c|}{Name} & \multicolumn{3}{c}{Lie brackets in addition to \eqref{eq:kla-generic} and \eqref{eq:bargmann-generic}}\\\toprule
    \widehat\g & [\B,H] = \P & & \\
    \widehat \n^+ = \b^+_{\gamma = -1} & [\B,H] = \P & [H,\P] = -\B & \\
    \widehat \n^- = \b^-_{\chi = 0}  & [\B,H] = \P & [H, \P] = \B &  \\\midrule
    \b^+_{\gamma \in (-1,1)} & [\B,H] = \P & [H,\P] = \gamma \B + (1+\gamma) \P & [H,Z] = (1+\gamma) Z\\
    \b^0 & [\B,H] = \P &  [H,\P] = \B + 2 \P & [H,Z] = 2 Z\\
    \b^-_{\chi>0} & [\B,H] = \P & [H, \P] = (1+\chi^2) \B + 2 \chi \P & [H,Z] = 2 \chi Z \\\bottomrule
  \end{tabular}
\end{table}

The Klein pairs $(\b,\h)$ are such that $\h \cong \iso(d)$ is the Lie
subalgebra spanned by $L_{ab},B_a$.  All Klein pairs are reductive,
with complementary subspace $\m$ spanned by $P_a, H, Z$.  Let $\pi^a,
\eta, \zeta$ be the canonical dual basis for $\m^*$.  All Klein pairs
$(\b,\h)$ share the same $\h$-invariant tensors in the linear isotropy
representation.  Up to scale they are given by the vector $Z \in \m$,
the covector $\eta \in \m^*$ and the lorentzian inner product
$h := \delta_{ab}\pi^a \pi^ b - 2 \eta \zeta \in \odot^2\m^*$.  Notice that
$Z$ is null relative to the inner product.  The vector field
corresponding to $Z$ is not only Killing but actually parallel with
respect to the Levi-Civita connection of the metric corresponding to
$h$.  The Klein pairs $(\b,\h)$ are geometrically realisable and
correspond to ($d+2$)-dimensional homogeneous lorentzian manifolds
and, by construction, their null reduction along $Z$ is $(\k,\h)$,
where $\k = \b/\RR Z$, abusing notation slightly and
denoting by $\h$ both the subalgebra of $\b$ and its isomorphic image
in $\k$.  Comparing with Table~\ref{tab:spacetimes}, it is evident
that the reduced Klein pairs $(\k,\h)$ are precisely the galilean
Klein pairs.  It should be remarked, however, that in the torsional
cases, the galilean structure obtained via the null reduction is not
the invariant one, since in those cases the null vector is not
invariant.

I conclude this section with an observation: the lorentzian metrics
on the geometric realisations of $(\widehat\g,\h)$ and $(\b^+_{\gamma
  = 0},\h)$ are flat and one can show that as Newton--Cartan
geometries, torsional de~Sitter--Galilei spacetime (for $\gamma = 0$)
and Galilei spacetime itself are isomorphic, although their
descriptions as kinematical Klein geometries are not.

\subsection{Carrollian geometry}
\label{sec:carrollian-geometry}

Carrollian geometry is in a sense dual to Newton--Cartan geometry.
Again let $V = \RR^{d+1}$ with elementary basis
$\be_0,\be_1,\dots,\be_d$ for $V$ and canonically dual basis
$\alpha^0,\alpha^1,\dots,\alpha^d$ for $V^*$.  Let $G \subset
\GL(V)$ denote the subgroup which leaves invariant $\be_0 \in V$ and
$\sum_{a=1}^d \alpha^a \alpha^a  \in \odot^2V^*$.  Explicitly,
\begin{equation}
  \label{eq:G-car}
  G = \left\{
    \begin{pmatrix}1 & \bv^T \\ \bzero & A\end{pmatrix} ~ \middle | ~  \bv
  \in \RR^d,~ A \in O(d)\right\} < \GL(d+1,\RR),
\end{equation}
with Lie algebra
\begin{equation}
  \label{eq:g-car}
  \g = \left\{
    \begin{pmatrix}0 & \bv^T \\ \bzero & A\end{pmatrix} ~ \middle | ~  \bv
    \in \RR^d,~ A \in \so(d)\right\} < \gl(d+1,\RR).
\end{equation}
The group of a carrollian structure is abstractly isomorphic to that
of a Newton--Cartan structure, both being isomorphic to the euclidean
group $\ISO(d) \cong \Ort(d) \ltimes \RR^d$.  Crucially, however, they
are not conjugate inside $\GL(V)$, so that they lead to different
geometries with different characteristic tensor fields.

The characteristic tensor fields of a carrollian geometry are a
nowhere-vanishing vector field $\xi \in \eX(M)$ (the
\textbf{carrollian vector field}) and a corank-one
positive-semidefinite symmetric tensor field $h  \in \Gamma(\odot^2
T*M)$ (the \textbf{spatial metric}), with $h(\xi,-) = 0$.  The triple
$(M,\xi,\h)$ defines a \textbf{(weak) carrollian geometry}, which can
be promoted to a \textbf{carrollian geometry} by the addition of an
adapted affine connection $\nabla$.  If the structure group
reduces to the connected component $G_0$ of G -- e.g., if $M$ is
simply connected -- then there is an addition a nowhere-vanishing top
form $\mu \in \Omega^{d+1}(M)$.

As shown in \cite{Figueroa-OFarrill:2020gpr}, the kernel and cokernel
of the Spencer differential are isomorphic as $G$-modules to
$\odot^2\Ann \be_0$, with $\Ann \be_0 \subset V^*$ the annihilator of
$\be_0$.  The isomorphism $\coker\d \cong \odot^2\Ann \be_0$ is
induced from a $G$-equivariant linear map $\Hom(\wedge^2V , V) \to
\odot^2 \Ann\be_0$ which induces in turn a bundle map $\Omega^2(M,TM)
\to \Gamma(\odot^2\Ann\xi)$ under which the torsion of an adapted
connection is mapped to $\eL_\xi h$, the Lie derivative of the spatial
metric along the carrollian vector field.

For $d>1$, there are four $G$-subbundles of $\odot^2\Ann\xi$ and hence
four classes of carrollian geometries depending on the intrinsic torsion of
adapted connections. For reasons which will only become clear after
discussing the geometric realisation of carrollian geometries as null
hypersurfaces in lorentzian manifolds, I shall refer to them as follows:
\begin{enumerate}
\item \textbf{totally geodesic}, if $\eL_\xi h = 0$;
\item \textbf{minimal}, if $\eL_\xi \mu = 0$;
\item \textbf{totally umbilical}, if $\eL_\xi h = f h$ for some $f \in
  C^\infty(M)$; and
\item \textbf{generic}, otherwise.
\end{enumerate}
If $d=1$ there are only two submodules and hence two carrollian
structures: either $\eL_\xi h = 0$ or not.

\subsubsection{Null  hypersurfaces}
\label{sec:null-hypersurfaces}

As shown in \cite{Duval:2014uoa, Hartong:2015xda}, a null hypersurface
in a lorentzian manifold admits a carrollian structure.  Let us review
these results here.  For more details on null hypersurfaces, see
\cite{MR886772,MR1777311}.

Let $(N,g)$ be a lorentzian manifold and $M \subset N$ an embedded
hypersurface such that the restriction of $g$ to $M$ is degenerate; in
other words, $M$ is a \textbf{null hypersurface} of $N$.  Because of
the lorentzian signature of the metric, the restriction of $g$ to $M$ must
be positive-semidefinite and of corank $1$.  This implies that there
exists a future-directed, nowhere-vanishing null vector
$\xi \in \Gamma(TM)$ such that for all $p \in M$,
$T_p M = \xi_p^\perp$.  Integral curves of $\xi$ can be reparametrised
in such a way that they are null geodesics.  They are called the
\textbf{null geodesic generators} of the null hypersurface $M$.  The
null vector $\xi$ defines a (trivial) line bundle $L \subset TM$,
which is independent of the choice of $\xi$: indeed, $L = TM^\perp
\cap TM$ without reference to $\xi$.  The metric $g$ defines a
riemannian structure $h$ on the quotient vector bundle $TM/L$ over $M$
via $h(\overline X, \overline Y) = g(X,Y)$ where $\overline X = X \mod
L$.  It is easy to see that $h$ is well-defined since $g(X,Y)$ only
depends on the residue classes of $X,Y$ modulo $L$.

With some abuse of notation, let $h$ stand for the restriction of $g$
to $M$ as well for the induced riemannian structure on $TM/L$.  The
triple $(M,\xi,h)$ is a (weak) carrollian structure.  The
classification of carrollian structures via their intrinsic torsion
corresponds to the analogue for null hypersurfaces of the
classification of hypersurfaces in riemannian geometry, as I now
explain.

The Levi-Civita connection of $(N,g)$ defines a \textbf{null
  Weingarten map} $b : TM/L \to TM/L$ defined by
\begin{equation}
  b(\overline X) = \overline{\nabla_X \xi}.
\end{equation}
It manifestly depends on the choice of $\xi$, but notice that if $f
\in C^\infty(M)$ is a positive smooth function so that $\widetilde\xi
:= f \xi$ is another generator of $L$, then $\nabla_X (f\xi) = f
\nabla_X \xi \mod L$.  In particular, the null Weingarten map at $p$
depends only on the value of $\xi$ at $p$.

Notice that if $X,Y \in \Gamma(TM)$, $[X,Y] \in \Gamma(TM)$ and hence
$g(\xi,[X,Y]) = 0$. Therefore,
\begin{equation}
  \begin{split}
      g(\nabla_X \xi, Y) - g(X,\nabla_Y\xi) &= X g(\xi, Y) -
      g(\xi,\nabla_X Y) - Y g(X,\xi) + g(\xi, \nabla_Y X)\\
      &= -g(\xi, \nabla_X Y - \nabla_Y X)\\
      &= -g(\xi, [X,Y]) = 0.
  \end{split}
\end{equation}
Hence the \textbf{null second fundamental form} $B(\overline X,
\overline Y) := h(b(\overline X), \overline Y) = g(\nabla_X\xi, Y)$ is
a well-defined symmetric form on $TM/L$.  By analogy with the theory
of hypersurfaces in riemannian geometry, the null hypersurface $M$ is
said to be
\begin{itemize}
\item \textbf{totally geodesic}, if $B=0$;
\item \textbf{minimal}, if $\tr b = 0$;
\item \textbf{totally umbilical}, if $B = f h$ for some $f \in
  C^\infty(M)$; and
\item \textbf{generic}, otherwise.
\end{itemize}
Notice that although $B$ depends on $\xi$, it does so via
multiplication by a positive function and hence the above conditions
are independent on the choice of $\xi$.

A short calculation shows that the restriction of $\eL_\xi g$ to $M$,
denoted $\eL_\xi h$, agrees (up to an inconsequential factor of
$\tfrac12$) with the null second fundamental form of the hypersurface,
showing that the above classification of null hypersurfaces
corresponds to the classification of carrollian structures via their
intrinsic torsion and explains the names given to the carrollian
structures in Section~\ref{sec:carrollian-geometry}.

\subsubsection{Bundles of scales of conformal structures}
\label{sec:bundles-scales}

Another natural source of carrollian geometries are bundles of scales of
conformal structures.  (See, e.g., \cite{Curry:2014yoa} for a recent
review.)  The fundamental example is the conformal sphere, whose
bundle of scales can be identified with the future (or past) deleted
lightcone, which is a null hypersurface in Minkowski spacetime.

Let $(S^{d-1},g)$ denote the round ($d-1$)-sphere, thought of as the
unit sphere in $\RR^d$ with the euclidean metric.  Let $[g]$ denote
the set of metrics on $S^{d-1}$ conformal to the round metric:
\begin{equation}
  [g] = \left\{ \Omega^2 g ~ \middle | ~ \Omega \in C^\infty(S^{d-1})\right\}.
\end{equation}
Pick a point $x \in S^{d-1}$.  The value $g_x \in \odot^2 T^*_xS^{d-1}$
at $x$ of the round metric, defines a ray $Q_x = \{ \lambda^2 g_x ~
\mid ~ \lambda \in \RR^+\} \subset \odot^2 T^*_x S^{d-1}$.  The
collection $ Q = \sqcup_{x \in S^{d-1}} Q_x$ of all such rays can be
given a differentiable structure making $\pi : Q \to S^{d-1}$ into a
smooth principal $\RR^+$-bundle, with $\sigma \in \RR^+$ acting as
$\lambda^2 g_x \mapsto \sigma^2\lambda^2 g_x$.  Let $G$ be the
connected component of $\SO(d-1,1)$.  Then $G$ acts transitively on
$S^{d-1}$ via conformal transformations and it therefore acts on $Q$.
This action is also transitive and the stabiliser of $g_x \in Q$ is
isomorphic to the euclidean group $\ISO(d-1)$.  In fact, $Q$ is
$G$-equivariantly diffeomorphic to the future deleted lightcone $\LL
\subset \RR^{d,1}$, and the diffeomorphism sends $\lambda^2 g_x
\in Q$ to $(\lambda,\lambda x) \in \LL$, where $x$ is a unit-norm
vector in $\RR^d$, giving rise to the following commutative triangle:
\begin{equation}
  \begin{tikzcd}
    Q \arrow[rd] \arrow[rr,"\cong"]  & & \LL \arrow[ld]\\
    & S^{d-1} &
  \end{tikzcd}
\end{equation}
which exhibits the diffeomorphism $Q \to \LL$ as a bundle
isomorphism in addition to as an isomorphism of homogeneous
$G$-spaces.  The carrollian structure on $Q$ corresponds to the
carrollian structure on $\LL$: the carrollian vector field is the
fundamental vector field of the $\RR^+$-action and the corank-one
degenerate metric is the pullback via the projection $\pi : Q \to
S^{d-1}$ of the round metric on $S^{d-1}$.

Now consider a riemannian conformal manifold $(N,[g])$, with $[g] =
\left\{ \Omega^2 g ~\middle |~ \Omega \in C^\infty(N)~\text{nowhere
    zero}\right\}$ the conformal class of $g$.  Let $p \in N$ and let
$M_p = \left\{\lambda^2 g_p \middle | \lambda \in \RR^+\right\}
\subset \odot^2 T^*_pN$ be the ray in $\odot^2 T^*_pN$ defined by
$g_p$.  Then $M = \sqcup_{p \in N} M_p$ is the total space of a
smooth principal $\RR^+$-bundle $\pi : M \to N$ called the \textbf{bundle of
  scales of the conformal manifold $N$}.  Let $h = \pi^* g$ and let
$\xi$ be the fundamental vector field of the free right $\RR^+$-action
on $M$.  Then $(M,\xi,h)$ is a (weak) carrollian geometry.  This
construction of carrollian structures has played a rôle in some recent
work \cite{Herfray:2020rvq, Herfray:2021xyp,
  Figueroa-OFarrill:2021sxz, Herfray:2021qmp}.

\section{Conclusions, omissions and outlook}
\label{sec:conclusion-outlook}

Despite the prominent rôle played by lorentzian geometry in both
general relativity and quantum field theory, it is not the only
possible geometrical description of space and time.  There are
phenomenological reasons for considering non-lorentzian geometries
(e.g., condensed matter physics, hydrodynamics,...), but they are also
the geometries relevant to a largely unexplored edge of the Bronstein
cube in Figure~\ref{fig:bronstein} which might provide a new approach
to constructing a quantum theory of gravity.

In this short review, I have tried to give a flavour of some of the
better studied non-lorentzian geometries: galilean (or Newton--Cartan)
and carrollian.  I reviewed the classification of kinematical Lie
algebras (with spatial isotropy) and their associated Klein
geometries.  I observed that despite the plethora of Klein geometries,
they belong to a small class of Cartan geometries: lorentzian,
riemannian, galilean and carrollian.  We concentrated on the galilean
and carrollian geometries and defined them as $G$-structures,
identified their characteristic tensors and refined the classification
according to their intrinsic torsions.

I have omitted several topics.  One topic I did not cover is that of
the automorphisms of these non-lorentzian geometries.  For example,
let $(M,\tau,\lambda)$ be a weak Newton--Cartan geometry.  Its
automorphism group is the subgroup of diffeomorphisms of $M$
preserving $\tau$ and $\lambda$.  Contrary to what happens in
lorentzian geometry, the automorphism group need not be
finite-dimensional.  For example, as shown in \cite{Duval:1993pe} and
revisited in \cite{Figueroa-OFarrill:2019sex}, the Lie algebra $\a$ of
infinitesimal automorphisms of the galilean Klein geometries is an
infinite-dimensional Lie algebra known as the \textbf{Coriolis
  algebra}.  It is a split extension
\begin{equation}
  \begin{tikzcd}
    0 \arrow[r] & C^\infty(\RR_t, \iso(d)) \arrow[r] & \a \arrow[r] &
    \RR D \arrow[r] & 0,
  \end{tikzcd}
\end{equation}
where $C^\infty(\RR_t,\iso(d))$ is the Lie algebra of smooth functions
from the real line (with parameter $t$) to the euclidean Lie algebra
under the pointwise Lie bracket on which $D$ acts as the derivation
$\frac{d}{dt}$.  Of course, ``strengthening'' the structure by the
addition of an adapted connection reduces the size of the Lie algebra
of infinitesimal automorphisms to a finite-dimensional Lie algebra.
This is the well-known fact (see, e.g., \cite{MR2532439}) that
the automorphism group of a Cartan geometry is finite-dimensional.
Something similar, but more interesting, happens with carrollian
structures.  The Lie algebra of infinitesimal (conformal)
automorphisms of weak carrollian structures can be
infinite-dimensional and are, in fact, intimately linked with the
asymptotic symmetries of asymptotically flat lorentzian manifolds, the
so-called BMS group \cite{Bondi:1962px,Sachs:1962zza}, as shown
originally in \cite{Duval:2014uva,Duval:2014lpa} and further discussed
in \cite{Figueroa-OFarrill:2019sex}.

Another omission is supersymmetry.  I have stayed here in the realm of
classical geometry, but of course there is a notion of non-lorentzian
supersymmetry and supergeometry. The results here are far from
complete.  There are classifications of certain four-dimensional
kinematical Lie superalgebras and their associated Klein
supergeometries
\cite{Figueroa-OFarrill:2019ucc,Grassie:2020dga,Grassie:2021zgc},
extending earlier work on contractions of the Poincaré and
anti~de~Sitter superalgebras and referred to in those papers.

As for future work, an obvious next step is the study of natural
conditions which can be imposed on the curvature of the Cartan
connection of a kinematical Cartan geometry.  Some of these conditions
could have a variational origin, just like Ricci-flatness in
lorentzian geometry arises as the Euler--Lagrange equation of the
Einstein--Hilbert action.  Closer in spirit to the approach outlined
in this review is the construction of Cartan geometries via the
``gauging procedure''.  This is the Physics version of the
construction of a Cartan geometry from local data (i.e., from an atlas
of Cartan gauges, in the language of \cite{MR1453120}).  Doing so for
the Cartan geometry modelled on Minkowski spacetime leads to the
Hilbert--Palatini action and results in Ricci-flatness or, more
generally, the Einstein condition (in the presence of a cosmological
constant).  As shown in \cite{Wise:2006sm}, doing so for the Cartan
geometry modelled on (anti)~de~Sitter spacetime leads to the
MacDowell--Mansouri \cite{MacDowell:1977jt} formulation of Einstein
gravity.  Work is in progress with Emil Have, Stefan Prohazka and
Jakob Salzer to ``gauge'' some of the four-dimensional kinematical
Klein geometries of interest.

Another extension of the work reviewed here is to study geodesic
motion on the non-Lorentzian geometry.  Work is in progress with Can
Görmez and Dieter Van den Bleeken studying dynamics on the galilean
Klein geometries by studying the geodesics of the invariant
connections.  One could also study dynamics on these geometries via
Souriau's method of coadjoint orbits \cite{MR1461545}.

The unitary representation theory of the kinematical Lie
groups is largely unexplored, with the notable exceptions of the
classic work of Wigner and Bargmann \cite{MR1503456,MR24827} for the
Poincaré group and of Lévy-Leblond \cite{MR154608} for the Galilei
group.  This is an important problem which could benefit from the
attention of representation theorists.

To conclude, we would not like to finish without mentioning other
geometries which are closely related to the kinematical geometries
treated here: not just aristotelian (as already mentioned), but also
conformal, Lifshitz, Bargmann,... which are finding applications in an
expanding set of research areas.


\providecommand{\href}[2]{#2}\begingroup\raggedright\endgroup

\end{document}